\documentclass[10pt]{amsart}
\usepackage{a4,amssymb,times,amsmath,cancel,mathrsfs}
\usepackage{graphicx,color,epsfig}

\usepackage[all]{xy}
\usepackage{tikz}
\usetikzlibrary{calc,3d,arrows}

\def\Id{\operatorname{id}}

\let\cal\mathcal

\def\cC{{\cal C}}

\def\cF{{\cal F}}
\def\cG{{\cal G}}
\def\cH{{\cal H}}
\def\cI{{\cal I}}
\def\cJ{{\cal J}}

\def\cL{{\cal L}}

\def\cP{{\cal P}}

\def\qpol{{\cal Q}}

\def\cS{{\cal S}}

\def\eps{{\epsilon}}

\let\blb\mathbb
\def\normalize{\mathtt {normalize}}
\def\CM{{\blb {CM}}}
\def\MM{{\blb M}}
\def\MMM{{\blb {MM}}}
\def\NCCR{{\blb {NCCR}}}

\def\cGG{{\blb G}}

\def \ZZ{{\blb Z}}
\def \T{{\blb T}}

\def \NN{{\blb N}}
\def \Rl{{\blb R}}

\newcommand{\se}[1]{\begin{equation*}\begin{split}#1\end{split}\end{equation*}}

\newcommand{\cone}[1]{[ #1 ]}
\newcommand{\wis}[1]{{\text{\em \usefont{OT1}{cmtt}{m}{n} #1}}}

\newcommand{\C}{\mathbb{C}}
\newcommand{\N}{\mathbb{N}}
\newcommand{\Z}{\mathbb{Z}}
\newcommand{\R}{\mathbb{R}}
\def\ccL{{\mathscr L}}
\def\ccE{{\mathscr E}}
\def\ccK{{\mathscr K}}

\newcommand{\Rplus}{\mathbb{R_+}}

\newcommand{\mut}{\mathsf{mut}}

\newcommand{\rank}{\textrm{rank}}

\newcommand{\Spec}{\ensuremath{\mathsf{Spec}}}

\newcommand{\<}{\langle}
\renewcommand{\>}{\rangle}

\newcommand{\End}{\ensuremath{\mathsf{End}}}

\newcommand{\Ext}{\mathsf{Ext}}

\newcommand{\sign}{\mathsf{sign \,}}

\newcommand{\spec}{\wis{spec}}

\newcommand{\Span}{\mathrm{Span}}

\newtheorem{lemma}{Lemma}[section]

\newtheorem{theorem}[lemma]{Theorem}
\newtheorem{corollary}[lemma]{Corollary}

\theoremstyle{definition}

\newtheorem{definition}[lemma]{Definition}

\newtheorem{question}[lemma]{Question}

{

}

\theoremstyle{remark}

\newtheorem{remark}[lemma]{Remark}

\newcommand{\rep}{\ensuremath{\mathsf{rep}}}

\newcommand{\Proj}{\ensuremath{\mathsf{Proj}}}

\newcommand{\Mat}{\mathsf{Mat}}

\newcommand{\Rep}{\mathsf{Rep}}

\newcommand{\Hom}{\textrm{Hom}}

\newcommand{\GL}{\ensuremath{\mathsf{GL}}}

\newcommand{\carr}{\mathtt R}

\newcommand{\PM}{\cP}
\renewcommand{\Rplus}{\R_+}

\title{Generating toric noncommutative crepant resolutions}

\author{Raf Bocklandt}
\address{Raf Bocklandt\\
School of Mathematics and Statistics\\
Herschel Building\\
Newcastle University\\
Newcastle upon Tyne\\
NE1 7RU\\
UK}
\email{raf.bocklandt@gmail.com}

\xyoption{all}

\begin{document}
\begin{abstract}
We present an algorithm that finds all toric noncommutative crepant resolutions of a given
toric 3-dimensional Gorenstein singularity. The algorithm embeds the quivers of these algebras inside a real 3-dimensional torus such 
that the relations are homotopy relations. One can project these embedded quivers down to a 2-dimensional torus
to obtain the corresponding dimer models. We discuss some examples and use the algorithm 
to show that all toric noncommutative crepant resolutions of a finite quotient of the conifold singularity
can be obtained by mutating one basic dimer model. 
We also discuss how this algorithm might be extended to higher dimensional singularities.
\end{abstract}
\maketitle

\section{Introduction}\label{intro}

In \cite{nccrep} Van den Bergh introduced the notion of a noncommutative crepant resolution of a singularity, which is
an algebra that can act as a substitute of an ordinary commutative crepant resolution. 
A noncommutative crepant resolution of a singularity with coordinate ring $R$ is a homologically homogeneous algebra 
of the form $A=\End_R(T)$ where $T$ is a a reflexive $R$-module. In the case of three-dimensional terminal Gorenstein
singularities, the derived category of representations of $A$ is equivalent to the derived category of coherent sheaves of a commutative crepant resolution.

The algebra $A$ can be seen as a path algebra  of a quiver $Q$ with relations. 
The vertices of this quiver correspond to the direct summands of $T$ and the arrows
to a basic set of maps between them. One can also define a dimension vector $\alpha$ which assigns to each vertex 
the rank of the corresponding summand.
With these data the singular variety $\spec R$ can be recovered as a moduli space parameterizing $\alpha$-dimensional semisimple
representations and in many cases a commutative resolution can be constructed by taking a moduli space parameterizing $\alpha$-dimensional 
stable representations for some stability condition.

If $R$ is a toric singularity, we would like our noncommutative resolution to carry an additional toric structure.
This can be done by asking that all summands of $T$ are have rank 1 and are graded. This ensures that the moduli spaces can be constructed as
toric quotients of toric varieties. If this is the case we call $\End(T)$ a toric noncommutative crepant resolution.

An interesting problem is to classify all possible toric noncommutative crepant resolutions of a given singularity.
In this paper we will describe an algorithm that does this for any toric three-dimensional Gorenstein singularity.
We then use a method by Craw and Quintero-Velez \cite{crawvelez}to embed the quivers of these
algebras inside a real 3-torus, such that the relations are precisely the homotopy relations.
If the singularity is Gorenstein the quiver can be projected to
a 2-torus to obtain a dimer model. 
This is a combinatorial gadget that was originally introduced in string theory \cite{HanKen,FranHan,HanHer}.
We illustrate the power of the algorithm by looking at some special examples: singularities from reflexive polygons and
abelian quotients of the conifold. In these cases we can prove that all the dimer models corresponding to such a singularity
are connected by mutations.

Finally we study the possible generalization of this algorithm to the case where the toric singularity 
is not Gorenstein or has higher dimension. 
  
\section{Acknowledgements}

I want to thank Nathan Broomhead, Ben Davison and Travis Schedler for the interesting discussions
we had about this project during their visits to Newcastle.
I also want to thank Lutz Hille for the stimulating research stay 
I had in Bielefeld.

\section{Preliminaries}
\subsection{Noncommutative crepant resolutions}

Let $R$ be the coordinate ring of a singular variety $X$.
A proper surjective map $\pi :\tilde X \to X$ is a \emph{resolution} if it induces an
isomorphism on the function fields and $\tilde X$ is a smooth variety.

There are many different ways to resolve a singularity and therefore one wants to impose 
extra conditions. A resolution is called \emph{crepant} if the pullback of the canonical divisor of $X$ under $\pi$ is
the canonical divisor of $\tilde X$. In particular if $X$ is a Gorenstein singularity (i.e. the canonical divisor is trivial)
then $\tilde X$ must be a Calabi Yau variety (i.e. a smooth variety with a trivial canonical bundle).

In 2 dimensions crepant resolutions are unique, 
in 3 dimensions a singularity can have more than one crepant resolution but they are closely related.
If $\tilde X_1\to X$ and $\tilde X_2\to X$ are two crepant resolutions then
Bridgeland \cite{bridgeland} proved their bounded derived categories of coherent sheaves are the equivalent.

In \cite{nccrep} Van den Bergh introduced special algebras which can act as a noncommutative
analogue of a crepant resolution. An algebra $A$ is a \emph{noncommutative crepant resolution} of $R$ if
it satisfies 2 conditions
\begin{itemize}
 \item $A\cong \End(T)$ where $T$ is a finitely generated reflexive $R$-module (reflexive means $\Hom_R(\Hom_R(T,R),R)\cong T$),
 \item $A$ is homologically homogeneous i.e. all simples have the same projective dimension.
\end{itemize}
He used this definition to show that in 3 dimensions these algebras behave like crepant resolutions: 
\begin{theorem}[Van den Bergh]\label{VdB}
For a terminal 3-dimensional Gorenstein singularity, the bounded derived category of finitely generated modules of a noncommutative crepant resolution 
is equivalent to the bounded derived category of coherent sheaves of a crepant resolution and vice versa. 
Noncommutative crepant resolutions exist if and only if a commutative crepant resolution exists. 
\end{theorem}
From now on we will abbreviate the term commutative crepant resolution by CCR and noncommutative crepant resolution by NCCR. 

\subsection{Maximal modification algebras}
A nice way to construct NCCRs was introduced by Iyama and Wemyss in \cite{wemyss}. 

For a positively graded commutative ring $R$ with $\C=R/R^+$ and Krull dimension $n$
we will call a module $T$ \emph{Cohen-Macaulay (abbreviated as CM)} if $\Ext_R^i(\C,T)=0$ for all $i<n$.

An algebra $A$ is called a \emph{modification algebra} if it is of the form $\End_R(T)$ with $T$ a reflexive module
and $A$ is Cohen-Macaulay as an $R$-module. $\End_R(T)$ is called a \emph{maximal modification algebra (abbreviated as MMA)} if for every
reflexive module $U$ such that $\End_R(T\oplus U)$ is Cohen-Macaulay then $U$ is a direct sum of summands of $T$. 
In other words one cannot add anything new to $T$ such that its endomorphism ring stays Cohen-Macaulay.
If $\End(T)$ is a (maximal) modification algebra we will call $T$ a (maximal) modification module. 

An NCCR is always an MMA, but the converse is not true: it might be the case that a singularity has an MMA that does not
have finite global dimension. But as soon as a modification algebra has finite global dimension then it is an NCCR (see Lemma 4.2 in \cite{nccrep}).

In three dimensions we have the following:
\begin{theorem}[Iyama-Wemyss \cite{wemyss}]\label{MMAisNCCR}
If $R$ is a 3-dimensional Gorenstein singularity for which there exists an NCCR, then the NCCRs are precisely the MMAs.
\end{theorem}

\subsection{Toric geometry}
Let $M$ be a group isomorphic to $\Z^n$ and let $N=\Hom(M,\Z)$, we use $\<,\>$ denote the pairing between $M$ and $N$ and
between $M\otimes \Rl$ and $N\otimes \Rl$.
Now let $v_1,\dots, v_k$ be vectors in $N$ such that the $\Rplus v_i$
form the rays of an integral polyhedral cone $\sigma \subset N\otimes \Rl$. 
We also assume that $\lambda v_i \in \N \iff \lambda \in \Z$, so the vectors
are irreducible. We write $\sigma= \cone{v_1, \dots, v_k}:= \Rplus v_1+\dots + \Rplus v_k$.

The \emph{dual cone} $\sigma^\vee$ is defined as 
\[
 \sigma^\vee := \{ x\in M\otimes\Rl | \forall u \in \sigma: \<x,u\>\ge 0\}.
\]
If we intersect this cone with $M$ we get a semigroup and we can use this semigroup
to construct a semigroup algebra $R_{\sigma} := \C[\sigma^\vee \cap M]$. This ring is positively graded
and normal. 
A cone $\tau= \cone{{v_{i_1}},\dots,{v_{i_l}}}$ is called a \emph{face} of the 
cone if and only if there is an $m\in \sigma^\vee$ such that $\tau = \{x \in \sigma| \<m,x\>=0\}$.
A face $\tau$ of $\sigma$ will give an embedding $\sigma^\vee\subset \tau^\vee$ and a inclusion 
$R_\sigma \subset R_\tau$ and one can see $R_\tau$ as a localization of $R_\sigma$.
This gives an embedding $\Spec R_{\tau} \subset \Spec R_{\sigma}$.

A \emph{fan} $\cF$ is a collection of integral polyhedral cones in $N\otimes \Rl$ such that 
\begin{itemize}
\item if $\tau$ is a face of $\sigma\in \cF$ then $\tau \in \cF$ and 
\item if $\sigma_1,\sigma_2 \in \cF$ then their intersection is a face of both $\sigma_1,\sigma_2$.
\end{itemize}
If $\cF$ is a fan we get a variety $X_{\cF}$ by taking the direct limit of the $\Spec R_{\tau}$ with respect
to the embedding maps $\Spec R_{\tau} \subset \Spec R_{\sigma}$.

To every cone $\sigma$ we can associate its fan $\cF_{\sigma}=\{\tau| \tau \text{ is a face of }\sigma\}$.
The variety $X_\sigma := X_{\cF_\sigma}$ is $\Spec R_\sigma$. $X_\sigma$ will be singular if and only if the $v_1,\dots, v_k$ do 
not form a basis for $N$.
To construct a resolution of $X_{\sigma}$ we need to subdivide the cone $\sigma$ in smaller cones that all correspond to smooth
varieties. The fan $\cF$ containing all these cones will then give us a variety $X_\cF$ which is a resolution
of $X_{\sigma}$.

A toric singularity $X_\sigma$ with $\sigma= \cone{v_1,\dots,v_k}$ 
will be called \emph{Gorenstein} if and only if there is an element $m$ in $M$ such that
$\<m,v_i\>=1$ for all $v_i$. This means that all $v_i$ lie in a common hyperplane and they
form the vertices of a convex $k-1$-dimensional polytope. 

If the dimension of the singularity is $3$ this polytope is a convex polygon in
a plane with integral coefficients. Subdivide the polygon into triangles with size $\frac 12$ and take the cones over these triangles to obtain a subdivision of $\sigma$. 
The size of the triangles implies that the cones all give smooth varieties, so
this subdivision gives rise to a resolution of $X_\sigma$.
Because the vectors that generate the subcones
all lie in the same plane the resolution is Calabi-Yau and hence it is a CCR.
  
\subsection{Quivers and Dimers}
As usual a \emph{quiver} $Q$ is an oriented graph. We denote the set
of vertices by $Q_0$, the set of arrows by $Q_1$ and the maps $h,t$
assign to each arrow its head and tail.
A \emph{nontrivial path} $p$ is a sequence of arrows $a_1\cdots a_k$
such that $t(a_i)=h(a_{i+1})$. A \emph{trivial path} is just a vertex.  
We will assume that our quivers are strongly connected i.e. there exist
a path from each vertex to each other vertex.

The \emph{path algebra} $\C Q$ is the complex vector space with as
basis the paths in $Q$ and the multiplication of two paths $p$, $q$ is
their concatenation $pq$ if $t(p)=h(q)$ or else $0$.
The span of all paths of nonzero length form an ideal which we denote by $\cJ$.
A \emph{path algebra with relations} $A=\C Q/\cI$ is the quotient of a
path algebra by a finitely generated ideal $\cI \subset \cJ^2$. 

We will call a path algebra with relations $\C Q/\cI$ \emph{positively
graded} if there exists a grading $\carr: Q_1 \to \Rl_{>0}$ such that
$\cI$ is generated by homogeneous relations.

A \emph{dimer model} $\qpol$ is a strongly connected quiver $Q$ enriched with a 2 disjoint sets of cycles of length at least $3$: $Q^+_2$ and $Q_2^-$, such that
\begin{itemize}
\item[DO] {\bf Orientability condition.} Every arrow is contained exactly once in one cycle in $Q_2^+$ and once in one in $Q_2^-$. 
\item[DM] {\bf Manifold condition.} The incidence graph of the cycles and arrows meeting a given vertex is connected.
\end{itemize}
Considering the cycles as polygons and gluing them together by the arrows we get a compact surface.
We speak of a dimer model on a torus if this surface is a torus.

A \emph{consistent $\carr$-charge} is a grading $\carr:Q_1 \to \R_{>0}$ such that
\begin{itemize}
\item[R1]
$\forall c \in Q_2: \sum_{a \in c}\carr_a =2$,
\item[R2] 
$\forall v \in Q_0: \sum_{h(a)=v}(1-\carr_a) + \sum_{t(a)=v}(1-\carr_a)=2$.
\end{itemize}

Given a dimer model we can construct its \emph{Jacobi algebra}
\[
 A_\qpol := \C Q/\<p-q| \exists a \in Q_1: pa \in Q_2^+,qa \in Q_2^-\>
\]

The main results concerning dimers and NCCRs

\begin{theorem}[Broomhead-Bocklandt]
The Jacobi algebra of a dimer model is an NCCR of its center if and only if it admits a consistent $\carr$-charge. 
In this case the center is a toric 3-dimensional Gorenstein singularity.
\end{theorem}
\begin{remark}
Broomhead proved this theorem for geometrically consistent $\carr$-charges (i.e. $\forall a\in Q_1:\carr_a<1$) in \cite{broomhead}.
The method of this proof was later extended to include all consistent $\carr$-charges in \cite{consistency}.
\end{remark}

\begin{theorem}[Gulotta-Ishii-Ueda]\label{atleastonedimer}
For every toric 3-dimensional Gorenstein singularity, there is a dimer model with a consistent $\carr$-charge
that has this singularity as the center of its Jacobi algebra.
\end{theorem}
\begin{remark}
The proofs of Gulotta \cite{gulotta} and Ishii and Ueda \cite{ishiidimer} are constructive: they provide an algorithm to construct at least one dimer model.
Varying the starting parameters of the algorithm one can obtain more than one dimer model for a
given singularity, but it is far from clear whether you will obtain all possible dimer models.
One of the aims in this paper is to construct a new algorithm that produces all dimer models for a given
singularity.
\end{remark}

Not every NCCR of a toric 3-dimensional Gorenstein singularity comes from a dimer model, 
but if we restrict to toric NCCRs, a natural restriction we will introduce in the next section,
one has the following:

\begin{theorem}[Bocklandt \cite{consistency,bocklandtqp}]\label{consistency}
Every toric NCCR of a toric 3-dimensional Gorenstein singularity is isomorphic to
the Jacobi algebra of a dimer model that admits a consistent $\carr$-charge. 
\end{theorem}

\section{Toric NCCRs and MMAs}

Now we can combine the concepts of the preliminary sections to obtain the notion of 
\emph{toric} NCCRs and MMAs. 

Let $A=\End_R(T)$ be an NCCR of $R=\C[X]$. 
The standard way to construct a commutative crepant resolution goes as follows.
First one writes $A=\C Q/\cI$ as the path algebra of a quiver with relations.
The vertices of this quiver will be in one-to-one correspondence with
the direct summands of $T$. Let $T_v$ be the direct summand corresponding to vertex $v$

To this quiver we can associate a dimension vector $\alpha : Q_0 \to \N: v \mapsto \rank T_v$, which allows us
to define the scheme
$\Rep_\alpha Q=\oplus_{a \in Q_1}\Mat_{\alpha_{h(a)}\times \alpha_{t(a)}}(\C)$, 
which parametrizes all $\alpha$-dimensional presentations of the path algebra $\C Q$.
This scheme contains a closed subscheme $\rep_\alpha A$ that parametrizes all $\C Q$-representations
that are also $A$-representations.  
On both schemes there is an action of $\GL_{\alpha}=\prod_{v \in Q_0}\GL_{\alpha_v}(\C)$ and the
quotient $\rep_\alpha A/\!\!/\GL_{\alpha}$ is isomorphic to $X=\Spec R$ (this is because $A$ is an $R$-order, see \cite{lebruynbook}).

If $X$ is a toric variety we want the above story to happen in the toric world. This means that the
main component of $\rep_\alpha A$ is a toric variety and the group $\GL_{\alpha}$ is a torus.
This is the same as asking that all $T_v$ are graded by $M$ (where $R\subset \C[M]$) and have rank $1$.

\begin{definition}
An NCCR or MMA is a \emph{toric NCCR} (\emph{toric MMA}) if it is of the form
$\End_R(T)$ where $T$ is a direct sum of $M$-graded rank 1 $R$-modules.
\end{definition}
The main questions we are interested in is the following:
\begin{itemize}
 \item []
\item Given a toric singularity can we find all toric NCCRs/MMAs?
\item If the singularity is 3-dimensional and Gorenstein can we find the dimer models corresponding to these NCCRs?
\end{itemize}

If $T=\oplus_v T_v$ is a direct sum of reflexive graded rank one
modules and $U$ is any reflexive graded rank one module then $\End(T\otimes_R U)\cong\End(T)$ (because $U$ is invertible), 
so by tensoring with $\Hom(T_v,R)$ we can make sure $R$ itself is one of the summands of $T$.
Now $\End(T)=\oplus_{v,w} \Hom(T_v,T_w)$ and $\Hom(R,T_v)=T_v$ so if $\End(T)$ is an MMA then it
is Cohen-Macaulay so all direct summands $T_v$ are also Cohen-Macaulay. From \cite{gubeladze} we have the following theorem 

\begin{theorem}[Gubeladze]\label{finite}
For a given toric singularity there are only a finite number of graded rank 1 Cohen-Macaulay modules
up to isomorphism.
\end{theorem}
\begin{corollary}
For any toric singularity there are only a finite number of toric MMAs and NCCRs. 
\end{corollary}

This finiteness result implies that the problem of classifying all toric NCCRs is feasible (at least in principle) and the first step
will be to obtain a list of all graded rank 1 Cohen-Macaulay modules.

\section{Graded rank 1 Cohen-Macaulay modules}
In this section we discuss a method to determine the graded rank 1 Cohen-Macaulay modules. The main ideas go back to Stanley, Bruns and Herzog in 
\cite{Stanley} and \cite{CMrings}. The method has also been described by Perling in more detail in \cite{perling}, \cite{perlingts}.

\subsection{Singular and local homology}
Let $S = \{v_1,\dots v_k\}$ a set of primitive vectors in $N$ that are the extremal rays of an $n$-dimensional cone $\sigma$.
For every $\varsigma \in \cF_\sigma$ we fix an ordered basis $B_\varsigma$ for the space $\{x \in M \otimes \R | \forall u \in \varsigma : \<u,x\>=0\}$ 
consisting of elements in $\sigma^\vee\cap M$.
These bases can be used to define an incidence function
\[
\eps(\varsigma_1,\varsigma_2) = \begin{cases}
                  0 & \text{if $\varsigma_2$ is not a codimension $1$ subspace of $\varsigma_1$}\\
      \sign \det G& \text{if $u \in B_{\varsigma_2} \setminus \Span B_{\varsigma_1}$ and $(G B_{\varsigma_1}, Gu)=B_{\varsigma_2}$.}
                  \end{cases}
\]
Let $\cG$ be any subset of $\cF_\sigma$ and define
$\cG_i$ to be the subset of $\cG$ containing all the faces of codimension $i$ ($|B_\sigma|=i$) and $\Z^{\cG_i}$ the free abelian group
generated by this set. Furthermore let
\[
d : \Z^{\cG_i} \mapsto \Z^{\cG_{i+1}} : \varsigma \mapsto \sum_{\kappa \in \cGG_{i+1}} \eps(\varsigma,\kappa) \kappa.
\]
It is easy to check that this turns
\[
\cS_\cG := \Z^{\cG_\bullet},d
\]
into a complex and the homology of this complex is called the \emph{singular homology of $\cGG$}.
\begin{lemma}
Let $\sigma$ be a cone.
The singular homology of the fan $\cF_\sigma$ is zero : $H(\cS_{\cF_\sigma})=0$.
\end{lemma}
\begin{proof} This is an easy consequence of the fact that a cone is contractible. \end{proof}

Let $R=\C[\sigma^\vee\cap M]$ be the ring defined by the cone $\sigma$, for any object $\varsigma$ in the fan $\cF_\sigma$ we put
\[
R_\varsigma := \C[\varsigma^\vee\cap M] = R[u_1^{-1},\dots ,u_k^{-1}]\text{ where $B_\varsigma=\{u_1,\dots, u_k\}$}
\]
We can also define a complex of $R$-modules $\cL_\bullet,d$ with
\[
\cL_i := \bigoplus_{\varsigma\in \cF_i} R_\varsigma
\]
for any $x= (x_\varsigma)_{\varsigma \in \cF_i}$ we set $dx=((dx)_\varsigma)_{\varsigma \in \cF_{i+1}}$ with
\[
(dx )_\kappa := \sum_{\varsigma \in \cF_i}\eps(\varsigma,\kappa) x_\varsigma
\]
where $x_\varsigma \subset R_\varsigma\subset R_\kappa$. 
This will give us again a complex and its homology is called the \emph{local homology of $R$}. If we tensor this complex with an $R$-module
$K$ and take the homology we get the \emph{local homology of $K$}.

\begin{theorem}\cite{CMrings}[Corollary 6.2.6]
An $M$-graded $R$-module $K$ is Cohen-Macaulay if and only if
the $i^{th}$ local homology of $K$ is trivial for $i=0,\dots,d-1$, where $d$ is the dimension of $R$.
\end{theorem}

We are now going to apply this to graded rank $1$ reflexive submodules of $R$. 

\subsection{Graded reflexive rank 1 modules}

Starting from a cone $\sigma=\cone{v_1,\dots,v_k}$,
we define for each $k$-vector of integers $b =(b_1,\dots, b_k) \in \Z^k$
\[
\T(b) = T_b :=\{x \in M| \<\alpha_i,x\> \ge b_i\} \text{ and }T(b)=\Span_\C \T(b). 
\]
The first can be considered as a semigroup module over the semigroup $\T(0)=\sigma^\vee \cap M$, while
the second is an algebra module over the algebra $R=T(0)$.
For different $b$'s the $T(b)$'s might be isomorphic as $R$-modules.
Indeed it is well known \cite{perling} that
\begin{lemma} 
\begin{itemize}
\item[]
 \item An $M$-graded $R$-submodule of $\C[M]$ is reflexive if and only if it is of the form
$T(b)$ for some $b \in \Z^k$.
 \item
If $b,b'\in \Z^k$ then $T(b)\cong T(c)$ as $R$-modules if and only if there is an $m \in M$ such that
$b_i = b_i' + \<m,v_i\>$.
\end{itemize}
\end{lemma}

Because all Cohen-Macaulay modules are reflexive, the problem of determining graded rank $1$ Cohen-Macaulay
modules reduces to checking which $T(b)$ are Cohen-Macaulay or in other words when 
$\cL_\bullet \otimes_R T(b)$ has trivial homology.

If $m\in M$ we can look at the degree $m$ part of the complex $\cL_\bullet \otimes_R K$.
\begin{lemma}[Bruns-Herzog \cite{CMrings}, Perling \cite{perling}]
For any module $T(b)$ we have that
$(\cL_\bullet \otimes_R T(b))_m  = \cS_{\cG_m}$ where 
\[
\cG_m = \{\varsigma \in \cF | \forall v_i \in \varsigma\cap S: \<m,v_i\>\ge b_i\}
\]
\end{lemma}
To check whether a given $T(b)$ is Cohen-Macaulay we must check that  $H((\cL_\bullet \otimes_R T(b))_m)=0$
or equivalently that $H(\cS_{\cG_m})=0$ for all the $\cG_m$ defined above.

For a given $m\in M$ and a $T(b_1,\dots,b_k)$ we denote the signature of $m$ for $v_i$ $+$ if
$\<m,v_i\>\ge b_i$ (here positive includes zero)  and $-$ otherwise.

Given $T(b_1,\dots, b_k)$ we define for every given signature $s \in \{+,-\}^k$  the \emph{cell} $C^s$ as
\[
C^s := \left\{ x \in M\otimes \R | \<v_i,x\> \begin{cases}
                                          \ge 0 &s_i=+\\
                                          <0 &s_i=-\\
                                         \end{cases}
\right\}
\]
The \emph{homology of the cell} is defined as the homology of $\cS_\cG$ where
$\cF^s$ contains only the faces in $\cF$ spanned by vectors with positive sign.
\[
\cF^s := \{\varsigma \in \cF| \forall i\in [1,k]: v_i\in \varsigma \implies s_i=+\} 
\]
              
With this definition $T(b_1,\dots,b_k)$ is Cohen-Macaulay if and only if
the cells with nontrivial homology in degrees $i=0,\dots,d-1$ contain no lattice points of $M$. 
In general checking the Cohen-Macaulay property reduces to checking whether integral solutions to
a large set of inequalities exist (see also \cite{perling}[Theorem 7.2]). This can become quite complex.
In dimension 3 things can be done to simplify the problem a lot.

\subsection{Some facts about Cohen-Macaulay modules in dimension $3$}\label{somefacts}
Let $R_\sigma$ be a toric singularity with cone $\sigma=\cone{v_1,\dots, v_k}\subset \Rl^3$.
We can choose a plane that intersects the cone transversally in
a polygon $\cP$. Without loss of generality we can choose the plane given by $z=1$ (i.e. the third coordinate is 1). 
If the singularity is Gorenstein
we can assume all vectors $v_1,\dots v_k$ lie in this plane. In general the vertices of the polygon $\cP$ are of the
form $\frac {v_i}{c_i}$ where $c_i \in \NN$ and (all $c_i=1$ in the Gorenstein case).

The standard orientation of $z=1$-plane will give us cyclic (clockwise) order to the vectors $v_i$, which we will identify with the index $i \in \Z/k\Z$.
Using this cyclic order, it is easy to see that
\[
\cF := \{ 0, \cone{v_1,\dots, v_k}\} \cup \{\cone{v_i}| i \in \Z/k\Z\} \cup \cup \{\cone{v_i,v_{i+1}}| i \in \Z/k\Z\}
\]
This has an interesting consequence:
\begin{lemma}\label{segment}
If $s \in \{+,-\}^k$ is a signature
then $\cF^s$ has trivial singular homology if and only if
the $i$ with positive signature $s_i$ consist of a sequence of numbers mod $k$.
\[
s_j=+ \iff j=i,i+1,\dots ,i+u \mod k
\]
\end{lemma}
\begin{proof}
It is clear that all $s$ of this form have trivial homology in degrees $0,1,2$.
If $s$ is not of this form $\cF^s$ consists of different connected components
so the zeroth homology will be nontrivial.
\end{proof}

\begin{theorem}\label{no-bounded-cells}
If the singularity is 3-dimensional a nonempty cell has nontrivial homology if and only if it is bounded.
\end{theorem}
\begin{proof}
Note that because no $3$ $v_i$ are sitting are linearly dependent, every $4$ planes $\<v_i,x\>=b_i$ will
either go through a common point or bound a tetrahedron.

Suppose a cell $C^s$ has nontrivial homology then the sequence of $v_i$ with positive signature
is not connected, so we can find $i_1<i_2<i_3<i_4$ with alternation signs $-+-+$.
Then $C^s$ is contained in
\[
B= \left\{x \in M\otimes \Rl| 
\<x,(-1)^{u}v_{i_u}\>\ge (-1)^u b_{i_u}, u=1,\dots,4
\right\}
\]
We will now show that if this set is nonempty it is a solid tetrahedron or a point.
So let $y$ be in this set, if $B$ is not a tetrahedron or a point the there would be a direction
$z$ such that $y+\lambda z \in B$ for all $\lambda\in \Rplus$. I.e. $\<z,(-1)^{u} v_{i_u}\>> 0$ for all
$u=1,\dots,4$. However because the cone spanned by $v_1$ and $v_3$ intersects the cone by $v_2$ and $v_4$, there are positive $\lambda_i$ such that
$\lambda_1 v_1 -\lambda_2 v_2 + \lambda_3 v_3 - \lambda_4 v_4=0$, taking the inner product with $z$ gives a contradiction. 

On the other hand if a cell $C^s$ has trivial homology, the cone $\sigma$ contains two cones: The one generated by the $v_i$ with positive signature,
the one generated by the ones with negative signature. Because of lemma \ref{segment} these 2 cones only intersect in the zero and therefore we can find a plane through the origin in $N\otimes \R$ 
such that the two cones lie on different sides of the plane.
The normal to this plane in the direction of the positive cone will give us an element $z \in M\otimes \Rl$ such that
$\<z,v_i\>>0$ for $v_i$ with positive signature and $\<z,v_i\><0$ for $v_i$ with negative signature.
Therefore if $y \in C^s$ then $y+\lambda z\in C^s$ for all $\lambda>0$ so $C^s$ is not bounded.
\end{proof}

The previous theorem can be used to show that the Cohen-Macaulay property remains true if we remove one of the $v_i$ from the cone of the singularity.
\begin{corollary}\label{remove-one}
If $T(b_1,\dots,b_k)$ is a Cohen-Macaulay module for the 3-dimensional singularity
generated by $v_1,\dots, v_k$, then $T(b_1,\dots,b_{i-1},b_{i+1},\dots,b_k)$
is a Cohen-Macaulay module for the singularity
generated by $v_1,\dots,v_{i-1},v_{i+1},\dots,v_k$.
\end{corollary}
\begin{proof}
Any cell $C^{\hat s}$ for the singularity
generated by $v_1,\dots,v_{i-1},v_{i+1},\dots,v_k$ is the union of two cells $C^s$ (one of which might be empty)
for the original singularity. If $C^{\hat s}$ was bounded and contained lattice points, then
both $C^s$ would be bounded and at least one of them will contain a lattice point.
\end{proof}

In the other direction we can show that if we add an extra ray to the singularity every Cohen-Macaulay module lifts to
a discrete interval of Cohen-Macaulay modules.
\begin{corollary}\label{interval}
Let $v_1,\dots v_k$ the vectors in cyclic order for a given 3-dimensional singularity and $k>3$.
If $T(b_1,\dots,b_{i-1},b_{i+1},\dots,b_k)$
is a Cohen-Macaulay module for the singularity
generated by $v_1,\dots,v_{i-1},v_{i+1},\dots,v_k$, then there are numbers $l,u \in \Z$ such that
the module
$T(b_1,\dots,b_{i-1},b_i,b_{i+1},\dots,b_k)$ is Cohen-Macaulay for the original singularity
if and only if  $l\le b_i \le u$.
\end{corollary}
\begin{proof}
For each of the $m$ in $M$ we can look at
the signs of $\<m,v_{i-1}\>-b_{i-1}$ and $\<m,v_{i+1}\>-b_{i+1}$.

In order to be Cohen-Macaulay $\<m,v_{i}\>-b_{i}$ must have the either the same sign $(+,-)$ as
$\<m,v_{i-1}\>-b_{i-1}$ or $\<m,v_{i+1}\>-b_{i+1}$. This gives us a condition of the form
$b_i \ge \<m,v_{i}\>$ if both have positive sign and $b_i \le \<m,v_{i}\>-1$ if both have negative sign.
The intersection of all these conditions is an interval. This interval must be bounded because by theorem \ref{finite} the number of 
Cohen-Macaulay graded rank 1 modules is
a finite set.
\end{proof}
It is important to note that theorem \ref{no-bounded-cells} and its 2 corollaries only hold in dimension 3. 
We will discuss this further in section \ref{higherdims}.

\section{Endomorphism rings of reflexives and embedded quivers}\label{embedded}

An \emph{embedded quiver} consists of a manifold $\MM$ and a quiver $Q$ such
that $Q_0 \subset \MM$ and every $a \in Q_1$ is a continuous map $a:
[0,1] \to \MM$ with $h(a)=a(1)$ and $t(a)=a(0)$. In this way every
path in the quiver corresponds to a path in $\MM$.
The homotopy algebra of an embedded quiver is the quotient of the path
algebra by the ideal generated by the expressions $p-q$ where $p$ and
$q$ are homotopic paths.
The \emph{universal cover} $\tilde Q$ of an embedded quiver $Q$ is by definition 
the pullback of the embedded quiver under the universal covering map $\tilde \MM\to \MM$.

In general if $T=T(b^1)\oplus \dots \oplus T(b^\ell)$ is a direct sum of nonisomorphic graded reflexive rank 1 modules of a toric ring $R$ 
then $\End_R T$ is isomorphic to an embedded quiver. This can be seen as follows:
choose for each $T(b_i)$ a point in $p_i \in M\otimes \R$ and make sure $p_i-p_j \not \in M$ for $i\ne j$.
If we shift $\T(b_i)$ to $\T(b_i)+m$ with $m\in M$ we get a new graded reflexive rank 1 module, which is isomorphic
to $T(b_i)$. We will assign to this module the point $p_i+m$. 

Every homogeneous homomorphism from $T(b^i) \to T(b^j)$ corresponds to shifting 
$T(b^i)$ by a vector $m\in M$ such that it lies in $T(b^j)$. 
Therefore it makes sense to identify a homomorphism $\phi$
with the path $a_\phi$ in $M\otimes \R$ from $p_i+m$ to $p_j$.
We will denote the vector from $p_i+m$ to $p_j$ by $a_{\phi}$.

Because $p_i$ and $p_i+m$ correspond to isomorphic modules, we can take the 
quotient of $M\otimes \R$ by $M$ to obtain an $n$-dimensional torus.
Now let $Q_0$ be the set of the images of the $p_i$ in this quotient.
Fix a set of homogeneous morphism that generate $\End_R T$ and
we take for $Q_1$ the images of the paths $a_{\phi}$ in $(M\otimes R)/M$.

It is easy to check that $\End_R T$ is isomorphic to the homotopy algebra
of this embedded quiver and its universal cover is the infinite quiver we started from.

Although the embedding of the quiver depends on the choice of the $p_i$, different
choices will give rise to homotopic embeddings (i.e. there is a one parameter family of embedded quivers connecting them).
The aim is to find embeddings that look nice, i.e. such that the arrows do not wind more around the torus than needed.

One way to do this is the following construction, which appears in work by Craw and Quintero-Velez \cite{crawvelez}. 
Given vectors $v_1,\dots, v_k$ that generate a cone in $N$ we have a map from $\phi:\Z^k \to N$
by mapping $(b_1,\dots, b_k)$ to $\sum_i b_iv_i$ and a dual map
$\phi^T: M \to \Z^k: m \mapsto (\<m,v_1\>,\dots \<m,v_k\>)$. Now $\phi\phi^T: M\to N$ will give us an embedding of $M$ in $N$ and
$\phi\phi^T\otimes \R: M\otimes \R\to N\otimes \R$  will be an isomorphism.

The map $\phi$ allows us to associate to each graded reflexive submodule $T(b_1,\dots, b_k)$ a point in $N$ and
the inverse $(\phi\phi^T\otimes \R)^{-1}$ embeds $N$ inside $M\otimes \R$, so we can assign to 
$T(b_1,\dots, b_k)$ the point
\[
 \kappa(b_1,\dots, b_k) :=(\phi\phi^T\otimes \R)^{-1}\phi(b_1,\dots, b_k).
\]
This assignment has the property that if you shift $\T(b_1,\dots, b_k)$ by $m\in M$, the
corresponding point will also shift by $m$, which is precisely what we want. 
After factoring out 
$M$ we get a map $\bar \kappa: \Z^k \to M\otimes \Rl/M$.
It is clear that if all $T(b^j)$ are mapped by to different points in the torus under this map, 
$\End(\oplus_j T(b^j))$ will be isomorphic to homotopy algebra of the embedded quiver we constructed. 

In dimension 3 this is indeed true because of the following lemma.
\begin{lemma}\label{only-one-on-zero}
If $b_1,\dots,b_k$ are such that $\sum b_i v_i=0$ then $T(b_1,\dots,b_k)$ cannot be Cohen-Macaulay
unless all the $b_i$ are zero.
\end{lemma}
\begin{proof}
If $T(b_1,\dots,b_k)$ were CM, then the $b_i\ge 0$ must form a segment, so
the cone spanned by the $v_i$ with nonnegative $b_i$ and the
cone spanned by the $v_i$ with negative $b_i$ only intersect in the top.
So $\sum_{b_i\ge 0}b_iv_i$ and $\sum_{b_i< 0}-b_iv_i$ can only be the same
if they are both zero.
\end{proof}
\begin{corollary}\label{alldifferent}
If  $\End(\bigoplus_j T(b^j)$ is a 3-dim NCCR, then every $T(b^j)$ will be mapped by $\bar\kappa$ to a different point in the torus.
\end{corollary}
\begin{proof}
If $T(b^j)$ and $T(c^j)$ we mapped to the same point then $T(b^j-c^j)$ would be mapped to the zero, but we know the only Cohen-Macaulay sitting on the zero
is the trivial one so this would mean that $T(b^j)\cong T(c^j)$.
\end{proof}

\section{The Algorithm for dimension 3}

The first thing we need is a procedure to check whether two reflexives $T(b)$ and $T(c)$ are isomorphic.
If we fix a basis for $M$, we can look at the corresponding half open unit cube in $D=[0,1)^3 \subset M\otimes \Rl$.
This cube is a fundamental domain for the quotient $M\otimes \Rl/M$, so
for every $T(b)$ there is a unique $m\in M$ such that $\kappa(b)+m$ lies inside $D$.
The module $T(b+\phi^T(m))$ will be the unique graded reflexive module isomorphic to $T(b)$ that
is mapped to a point inside the unit cube. 
Therefore it make sense to define the function 
\[
 \normalize: \Z^k\to \Z^k: b \mapsto b+\phi^T(m).
\]
Two vectors in $\Z^k$ will define
isomorphic modules if and only if there image under the map $\normalize$ is the same.

The algorithm now consists of 3 steps
\begin{itemize}
\item[Step 1] Generate the set $\CM:= \{\normalize(b)| T(b) \text{ is Cohen-Macaulay}\}$.  
\item[Step 2] Find the maximal modifying sets: these the are subsets $S \subset \CM$ such that $0 \in S$ and $\forall b,c \in S, \normalize(b-c) \in \CM$.
\item[Step 3] For every maximal subset $S$ we construct the quiver of $\End(\oplus_{b \in S} T(b))$.
\end{itemize}

We will now discuss these 3 steps in a bit more detail.

\subsection{Generating the Cohen-Macaulays}

Because the 3-dimensional situation is special in the sense that the Cohen-Macaulay Property is maintained after removing one of the $v_i$, we can use an inductive procedure
to generate $\CM$. We order the $v_i$ cyclicly. The singularity generated by the first $3$ vectors
is a quotient singularity $\C^3/G$ where $G\subset \GL_3(\C)$ is abelian. This is
because the cone is simplicial. This implies that every reflexive is Cohen-Macaulay, so $\CM(v_1,\dots,v_3) = \normalize(\Z^3)$.

Given $\CM(v_1,\dots,v_i)$ we can construct $\CM(v_1,\dots,v_i,v_{i+1})$ as follows. 
For each $(b_1,\dots, b_i)$ we know from Lemma \ref{interval}
that we there is an interval $[l,u]$ such that the module $T(b_1,\dots, b_i,b_{i+1})$ is Cohen-Macaulay if and only if
$b_{i+1} \in [l,u]$. To find this interval we start with a given $z\in \Z$ check whether $(b_1,\dots, b_i,z)$ is Cohen-Macaulay or not.

This check is done by making looking at all $4$-tuples of planes $\<x,v_{i_j}\>-b_{i_j}-\eps=0$ where one of $i_j=i+1$ and check whether the tetrahedron they bound does not
contain any lattice points. (The $\eps$ is just a tiny number $(\eps<1)$ to compensate the fact that the cells are half open.)

\begin{itemize}
 \item[A]
If it is CM we put $\normalize(b_1,\dots, b_i,z)$ in $\CM(v_1,\dots,v_i,v_{i+1})$ and we check 
\[
(b_1,\dots, b_i,z+1),(b_1,\dots, b_i,z+2),\dots 
\]
until we get a vector that does not give a Cohen-Macaulay.
Then we do the same with 
\[
(b_1,\dots, b_i,z-1), (b_1,\dots, b_i,z-2), \dots 
\]
and by lemma \ref{interval} we can be sure all Cohen-Macaulays are found once we
get to a vector that is not Cohen-Macaulay.
\item[B]
If $(b_1,\dots, b_i,z)$ is not Cohen-Macaulay, there were lattice points in some of the tetrahedra.
Chose such a lattice point $m$ and change $z$ to $z'=\<m,v_{i+1}\>$ if $\<m,v_{i+1}\><z$ or $z'=\<m,v_{i+1}\>-1$ $\<m,v_{i+1}\>\ge z$ (to remove $m$ from the tetrahedron).
Now check whether $(b_1,\dots, b_i,z')$ is Cohen-Macaulay. 

\begin{itemize}
\item[B.1] If it is CM proceed according to paragraph A.
\item[B.2] If it is not CM we take again a lattice point in a tetrahedron an use it to modify $z'$ to $z''$ as above, and continue like this until we hit
a Cohen-Macaulay. 

If the sequence $z,z',z'',\dots$ is not monotone we stop. This means we cannot make a Cohen-Macaulay of the form $(b_1,\dots, b_i,z)$.
Note that the sequence can not continue monotonely indefinitely. Indeed, we start with a finite number of lattice points in the tetrahedra and changing $z$ will remove
some lattice points from the tetrahedra and maybe add some new ones, but in order to remove the new ones we have to move $z$ in the opposite direction.
\end{itemize}
\end{itemize}

\subsection{Finding the maximal subsets}

Choose any indexing $e_0,e_1,\dots$ of the elements in $\CM$ such that $e_0=0$.
We are going to construct generations of pairs of subsets of $\CM$, $(S,T)$ satisfying the following conditions:
\[
S\cap T=\{\} \text{ and } \forall b \in S:\forall b \in T, b-c,c-b \in \CM
\]
The first generation only contains the pair $(\{e_0\},\{e_1,e_2,\dots\})$. Given a generation we construct the next generation
by constructing for every pair $(S,T)$ in this generation and for every $e_j\in T$ with $j$ bigger than all indices in $S$, a new pair
\[
(S\cup \{e_j\}, \{c \in T\setminus\{e_j\}| e_j-c,c-e_j \in \CM \}).
\]
As $S$ increases and $T$ decreases, there will be a last nonempty generation.

The algorithm returns the set $\MMM$ containing all first entries of the pairs of this last generation.
\subsection{Constructing the quiver}\label{consq}

Following section \ref{embedded}, we construct an embedded quiver for every $S \in \MMM$. As vertices it will have
the set $Q_0 = \{\bar \kappa(b)| b \in S\}\subset M\otimes \Rl/M$ 
(which by corollary \ref{alldifferent} are all different) 
and the arrows $Q_1$ come from a minimal set of graded algebra generators for $\End_R(\oplus_{b \in S}T(b))$.
If the arrow $a$ corresponds to a graded homomorphism $T_b\to T_c$ we define $\vec a = \bar \kappa (c-b)$.
Using this notation we identify $a$ with the map 
\[
 a : [0,1] \to M\otimes \Rl/M : t \mapsto \kappa(b)+ t \vec a. 
\]

Many maximal modifying sets $S$ will give isomorphic homotopy algebras. In order to remedy this, 
we introduce the notion of affine equivalence.

An affine transformation $\Psi$ of $M\otimes \R$ (considered as an affine space)
is compatible with the quotient $M\otimes \R/M$ if it maps fibers to fibers.
This allows us to see it also as a map $\Psi : M\otimes \R \to M\otimes \R$.
We will call two embedded quivers $Q$ and $Q'$ in $M\otimes \R/M$ affine equivalent if there is an affine transformation
$\Psi$ compatible with the quotient 
such that
\[
a \in Q_1 \iff  \Psi a \in Q'_1
\]
For each $S$ the algorithm will now construct the quiver and check whether it is affine equivalent to one of the quivers already constructed.
If not this quiver will be added to the list $\NCCR$.
The complexity of checking affine equivalence of 2 embedded quivers is proportional to the number of arrows in the quiver because
an affine transformation is fixed once we know the image of one arrows.

\subsection{Effectiveness of the algorithm}
\begin{theorem}\label{effective}
Given the toric data $\{v_1,\dots,v_k\}\subset \Z^2\times\{1\}$ for a toric 3-dimensional Gorenstein ring $R$,
the algorithm stops and gives 3 lists $\CM$, $\MMM$ and $\NCCR$.
\begin{enumerate}
\item 
The list $\CM$ contains vectors $b \in \Z^k$ corresponding to all graded rank-1 Cohen-Macaulay modules $T(b_1,\dots,b_k)$ up to isomorphism.
\item
The list $\MMM$ will correspond to all maximal modifying modules that decompose as a direct sum of graded rank-1 Cohen-Macaulays, one of which is $R$ itself.
\item
The list $\NCCR$ contains all embedded quivers (inside the three-dimensional torus) for the toric NCCRs up to affine equivalence.
\end{enumerate}
\end{theorem}
\begin{proof}
\begin{enumerate}
 \item[]
\item 
Lemmas \ref{remove-one} and \ref{interval} and the monotony requirement for the sequence $z,z',\dots$, ensure that 
step 1 stops and $\CM$ contains all graded rank-1 Cohen-Macaulay modules.
Steps 2 and 3 also clearly stop.
 \item
If $R$ is Gorenstein, then by theorem \ref{atleastonedimer} we know that there is at least one toric NCCR, which by 
theorem \ref{MMAisNCCR} is given by a maximal modifying module (consisting of graded rank-1 Cohen-Macaulays because it is toric).
This maximal modifying module will be contained in $\MMM$. By construction the $\MMM$ contains sets of the same size, these must all be maximal modifying: 
if they were not, 
they were contained in a maximal modifying with more summands.
The corresponding NCCR would then have a higher rank Grothendieck group than the one toric NCCR we already had. 
This is impossible because by \cite{nccrep} all NCCRs are derived equivalent.
\item
This is by construction.
\end{enumerate}
\end{proof}

We can also recover the dimer models from the embedded quivers in the following way:
\begin{theorem}[Craw-Quintero-Velez \cite{crawvelez} Theorem 5.9]
Given any embedded quiver $Q$ in the list $\NCCR$, 
one can project it down to the twodimensional torus by forgetting the third coordinate.
After this projection the arrows will cut the 2-torus in polygons which are bounded by cycles.
These turn $Q$ into a dimer model.
\end{theorem}

As we already mentioned these dimer models are all consistent so this means there must exist consistent $\carr$-charges for
these dimer models.

\begin{theorem}\label{charges}
Choose a basis for $M$ such that the Gorenstein vector is $(0,0,1)$. 
Let $Q$ be any embedded quiver in the list $\NCCR$. 
For any $x$ in the interior of $\sigma^\vee$ we can construct a consistent $\carr$-charge for $Q$ by
putting $$\carr_a = 2\<x,\vec a\>/\<x, (0,0,1)\>.$$ 
where the inner product on $M\otimes \R$ comes from the basis and $\vec a$ is as defined in \ref{consq}.
\end{theorem}
The proof of this theorem is postponed to the end of this section because 
it uses the techniques of perfect matchings and zigzag paths.

A \emph{perfect matching} is a set of arrows $\PM \subset Q_1$ meeting every cycle
in $Q_2$ in precisely one arrow. It gives us a degree function on $A_\qpol$ by giving the arrows
in the matching degree $1$ and the others degree $0$. We denote the degree of an element $u\in A_\qpol$
under $\PM$ by $\PM(u)$. $\PM$ gives a degree function on $A_\qpol$ so by restriction also
on $Z(A_\qpol)=R_{\sigma}$. This degree function on $R_{\sigma}\subset \C[M]$ comes from an element in $N$
which we denote by $\bar \PM$. Note that different perfect matchings can have the same vector in $N$.

The convex hull of all these vectors $\bar \PM$ forms a polygon in $N$. 
We call a perfect matching extremal if its vector lies on a corner of this polygon.

\begin{theorem}[Broomhead]\cite{broomhead}
Let $\qpol$ be a consistent dimer model and $A_\qpol$ its Jacobi algebra.
There is a one to one correspondence between the corners of the polygon that defines
$R=Z(A_\qpol)$ and the extremal perfect matchings.
More precisely
\[
 R=R_{\sigma} \text{ with } \sigma=\cone{\bar \PM| \PM \text{ is an extremal perfect matching}}.
\]
\end{theorem}
\begin{remark}
The proof of the statement is done for geometrically consistent dimer models, but using section 8 of \cite{consistency} it works for
all consistent dimer models.
\end{remark}

We order the extremal prefect matchings cyclicly according to the polygon
$\PM_1,\dots,\PM_k$ and $v_i$ is the coordinate $(\PM_i(x),\PM_i(y),1)$ of the ith vertex of the polygon.
\begin{theorem}
Let $\qpol$ be a consistent dimer model and $A_\qpol$ its Jacobi algebra.
Fix a vertex $v \in \qpol_0$.
\[ A_\qpol= \End(\oplus_{u \in \qpol}T_u)\]
with 
$T_u := T(\PM_1(p_u),\dots, \PM_k(p_u))$ where $p_u$ is any path from $v$ to $u$.
\end{theorem}
\begin{proof}
Because $A_\qpol$ is a toric NCCR we know it is of the form $ A_\qpol= \End(\oplus_{u \in \qpol}T_u)$ for some $T_u$.
It is clear that $u A_\qpol v= \End_R(T_{v},T_{u})$. 

Suppose we put $T_v=R$.
Given any path $p_u$ from $v$ to $u$, we can embed $T_u=u A_\qpol v$ in $R=v A_\qpol v$ by multiplication
with $p_u$. This shows that $T_u \subset T(\PM_1(p_u),\dots, \PM_k(p_u))$.
To show that this is an equality we need to prove that for every extremal perfect matching $\PM_i$ there
is a path in $u A_\qpol v$ with $\PM$-degree zero. The proof of this follows from an adaptation of the proof of proposition 6.2
in \cite{broomhead} and theorem 8.7 in \cite{consistency} which states that 
given a homology class of paths from $u\to v$ (on the $2$-torus) we can find an extremal perfect matching $\PM_i$ 
and a path with that homotopy class $p\in vA_\qpol u$ such that $\PM_i(p)=0$. Now from the construction in these proofs it is clear
that every extremal perfect matching $\PM_i$ will occur if one varies the homology class.
This shows that $T(\PM_1(p_u),\dots, \PM_k(p_u))=T_u$.
\end{proof}
\begin{corollary}\label{arrows}
Let $\qpol$ be a consistent dimer model then
every arrow $a$ is contained in at least $1$ extremal perfect matching and at most $k-2$. Moreover
\[
\vec a = (\phi\phi^T\otimes \R)^{-1}\left(\sum_{i=1}^k \PM_i(a) v_i\right)  
\]
\end{corollary}
\begin{proof}
because $a$ is a path from $h(a)$ to $t(a)$ we know that 
if $T_{t(a)}=T(b_1,\dots,b_k)$ then
$T_{h(a)}=T(b_1+\PM_1(a),\dots, b_k+\PM_k(a))$.
Therefore $\vec a = (\phi\phi^T\otimes \R)^{-1}\left(\sum_{i=1}^k \PM_i(a) v_i\right)$.
Because $T_{t(a)}\ne T_{h(a)}$, $a$ is contained in at least one extremal perfect matching.

Every arrow $a$, is contained in a cycle $c$ in $Q_2$ of length at least $3$, because $\PM_i(c)=1$
for all $i$ and the two other arrows in $c$ are each contained in at least $1$ extremal perfect matching,
$a$ is contained in at most $k-2$ extremal perfect matchings.
\end{proof}

\begin{lemma}\label{conspm}
Let $\qpol$ be a consistent dimer model.
For every extremal perfect matching $\PM_i$ and every vertex we have the following property
\[
\sum_{h(a)=v} \PM_i(a) + \sum_{t(a)=v} \PM_i(a) = \#\{a\in Q_1| h(a)=v\} -1
\]
\end{lemma}
\begin{proof}
See the proof of theorem 8.7 in \cite{consistency}. 
\end{proof}

\begin{proof}[Proof of theorem \ref{charges}]
We can rewrite the equation in the previous lemma as 
\[
\sum_{h(a)=v} (\frac 12 -\PM_i(a)) + \sum_{t(a)=v} (\frac 12 - \PM_i(a)) = 1.
\]
If we multiply this with $v_i$, take the sum over $i$ and apply $(\phi\phi^T\otimes \R)^{-1}$ we get
\[
(\phi\phi^T\otimes \R)^{-1}\sum_i(\sum_{h(a)=v} (\frac 12v_i -\PM_i(a)v_i) + \sum_{t(a)=v} (\frac 12 - \PM_i(a)v_i)) = (\phi\phi^T\otimes \R)^{-1}\sum_i v_i.
\]
The special choice of basis gives
\[
\sum_{h(a)=v} ((0,0,\frac 12) -\vec a) + \sum_{t(a)=v} ((0,0,\frac 12) -\vec a) = (0,0,1).
\]
If we now take the inner product of this with $x$ and rescale by $2/\<x,(0,0,1)\>$, we
get the consistency condition for the $\carr$-charge. Note that $x$ needs to sit in $\sigma^\vee$
for the charges to be positive.
\end{proof}

\section{Examples}

\subsection{Reflexive polygons}
If a convex integral polygon in $\Z^2$ has exactly one internal lattice point 
it is called a reflexive polygon.
Up to an integral affine transformations, there are precisely $16$ reflexive polygons. 
We will choose the internal lattice point to be $(0,0)$ and let $v_1,\dots,v_k$ be the lattice points on the boundary of the polygon in cyclic order.
The fan 
\[
\{0, \cone{v_1}\,\dots, \cone{v_k}, \cone{v_1,v_2},\dots, \cone{v_k,v_1} \} 
\]
in $\R^2$ will define a projective smooth toric variety. This surface is a Fano surface 
if all $v_i$ are corners of the polygon and a weak Fano surface otherwise. 
To the fan we can associate a sequence of numbers
\[
 (a_1,\dots, a_k) \text{ such that } v_{i-1}+a_iv_i+v_{i+1}=0.
\]
And this sequence (up to cyclic shifts and inversion of the order) determines the isomorphism class of the Fano variety. 
Every divisor $v_i$ defines a line bundle $\ccE_i$ and these generate the Picard group. 
On the Picard group we have the intersection form:
\[
 \<\ccE_i,\ccE_j\> = \begin{cases}
1 &i=j\pm 1\\
a_i&i=j\\
0& |i-j|>1.
\end{cases}
\]

In \cite{perlinghille} Hille and Perling studied
full cyclic strongly exceptional sequences of line bundles associated to weak Fano surfaces.
Full strongly exceptional sequences are infinite sequences of line bundles $\dots, \ccL_i, \ccL_{i+1},\dots$ such that 
\begin{itemize}
 \item $\Ext^r(\ccL_i,\ccL_j)=\Ext^r(\ccL_j,\ccL_i)=0$ if $r>0$ and $i\le j< i+k$,
 \item $\Hom(\ccL_i,\ccL_j)=0$ if $i>j$,
 \item $\ccL_{i+k} = \ccL_{i}\otimes \ccK^{-1}$.
\end{itemize}
Here $\ccK$ is the canonical bundle and $k$ is the rank of the Grothendieck group (which in the toric case equals
the number of on dimensional cones in the fan).

\begin{theorem}[Hille, Perling \cite{perlinghille}]
Given a cyclic full strongly exceptional sequence $(\ccL_i)$ on a toric surface,
The sequence of numbers
\[
\<\ccL_{i+1}-\ccL_{i},\ccL_{i+1}-\ccL_{i}\>,\dots,\<\ccL_{i+k}-\ccL_{i+k-1},\ccL_{i+k}-\ccL_{i+k-1}\>
\]
corresponds to the sequence a new reflexive polygon.
\end{theorem}
We will call this sequence the type of the exceptional collection.
Note that this type can be different for different exceptional collections on the same Fano.
Hille and Perling also construct a table of types that can occur for each Fano surface.
 
This problem is closely related to our problem. Starting from a Fano variety 
we can also make a 3-dimensional Gorenstein singularity by embedding the polygon in $\R^3$ in the plane $z=1$. 
We denote its cone by $\sigma := \cone{\tilde v_1,\dots,\tilde v_k}$ where $\tilde v_i=(v_i,1)$.

The vector $(0,0,1) \in N$ gives rise to a grading on $\C[M]$ and hence also on $R_{\sigma}$ and
every reflexive module $T_b$. For this grading $\Proj R_{\sigma}$ gives us the Fano variety and
every graded module $T_b$ corresponds to a line bundle $\ccL_b$ over $\Proj R_{\sigma}$. 
Two isomorphic modules $T_b$, $T_c$ can give the non-isomorphic line bundles if their grading is different.
More specifically we have that
\[
 \ccL_b \cong  \ccL_c \iff T_b\cong T_c \text{ and } \<\kappa(b),(0,0,1)\>>\<\kappa(c),(0,0,1)\>.
\]

Let $T=T_{b_1}\oplus \dots \oplus T_{b_k}$ be a maximal modifying module for $R_{\sigma}$.
Without loss of generality we can assume that the $b_i$ are normalized ($\kappa(b_i)\in [0,1)^3$) and
ordered by increasing $\<\kappa(b_i),(0,0,1)\>$.
We now define a sequence of line bundles $(\cL_{j})_{j\in \ZZ}$: 
\[
 \ccL_{j+\ell k} := \ccL_{b_{j}+\ell (1,\dots, 1)} \text{ for $1\ge j\ge k$ and $\ell \in \ZZ$}.
\]

\begin{theorem}
This sequence is a cyclic full strongly exceptional sequence of line bundles. 
\end{theorem}
\begin{proof}
We can cover the Fano variety by its basic affine open sets coming from the fan. This cover
can be used to calculate the \v{C}ech complex $\check{\cC}^\bullet(\ccL_b)$ of a line bundle $\ccL_b$.
The $i^{th}$ component of this complex can easily be seen to be the 
$(0,0,1)$-degree zero part of the $(i+1)^{th}$ component of
$\cL\otimes T_b$. 
The fact that $T_{b-c}$ is CM for all $\ccL_b,\ccL_c \in \cH$ gives us that
$\Ext^i(\ccL_b,\ccL_c)=0$ if $i\ne 0,2$.

Now $\Hom(\ccL_b,\ccL_c)= [T_{c-b}]_0$ which is zero if $\<\kappa(c-b),(0,0,1)\>>0$ because
all homomorphisms have positive $\carr$-charge. 
Finally because of Poincar\'e duality $\Ext^2(\ccL_b,\ccL_c)=\Hom(\ccL_{c+(1,\dots,1)},\cL_b)^*$
and this is zero if $\<\kappa(b-(1,\dots,1)-c),(0,0,1)\>>0$ or
$\<\kappa(c-b),(0,0,1)\><1$.
\end{proof}
\begin{remark}
In the weak Fano case a similar result holds but one needs to tweak the method to extract $\cL_{b_i}$ from $T_{b_i}$:
because $\Proj R_\sigma$ is now singular and one needs to use its minimal resolution to pull back to the weak Fano.    
This can be done using a moduli construction. For more info see \cite{Ishii} and \cite{mozgovoy}.
\end{remark}

Using our algorithm we can determine all maximal modification modules of $R_\sigma$.
For each maximal modification module we choose a compatible index and let $w_i$ be the vertex in corresponding
the NCCR according to $T(b_i)$. From \cite{perlinghille} we get 
\[
 \<\ccL_i-\ccL_{i-1},\ccL_i-\ccL_{i-1}\>= \#\{\text{arrows from $w_{i-1}$ to $w_{i}$ in the NCCR}\}-2.
\]
The corresponding dimer models and their types can be found in the appendix and 
on {\tt www.algebra.ua.ac.be/dimers}.

\subsection{Mutations}\label{mut}
An important thing that has been noticed is that it is possible to turn a consistent dimer model $\qpol$
into another consistent dimer model $\mut_v\qpol$ such that $A_{\qpol}$ and $A_{\mut_v\qpol}$.
This process is called mutation. It orginates from cluster algebras \cite{fominzelev} and is applied to
algebras with a superpotential in \cite{derksenweyman}.

For dimer models on a torus the procedure restricts to the following construction.
Let $v$ be a vertex without loops or cycles of length $2$ and with exactly $2$ arrows $a_1,a_2$ leaving and two arrows $b_1,b_2$ arriving \footnote{the requirement on the amount of arrows arriving and leaving is not needed to define a mutation, however the mutated quiver will usually not be a dimer model any more}.
The toricly mutated dimer $\mut_v\qpol$ has the same vertices as $\qpol$ and 
to obtain a list of all arrows and faces for $\mut_v\tilde\qpol$ we apply the following procedures.
\begin{itemize}
 \item Replace the arrows $a_1,a_2,b_1,b_2$ by $a'_1,a'_2,b'_1,b'_2$ in the opposite direction.
 \item For each the path $a_ib_j$ add an arrow $u_{ij}$ with the same head and tail.
 \item If $a_ib_js$ is a face in $\qpol_2^{\pm}$ change it to $u_{ij}s$ and add a cycle
$b'_ja'_iu_{ij}$ to $\qpol_2^{\mp}$.
 \item If there are faces of length 2, remove them and stick the 2 faced that bound them together.
\end{itemize}
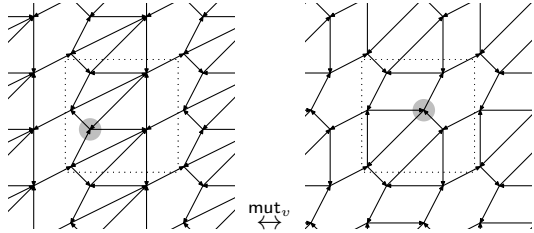
\begin{center}
\resizebox{3cm}{!}{\includegraphics{mut1.mps}}
$~~~\stackrel{\mut_v}{\leftrightarrow}~~~$
\resizebox{3cm}{!}{\includegraphics{mut2.mps}}
\end{center}

Another way to describe this toric mutation procedure, which relates to the viewpoint of the dimer as an MMA can be found
in \cite{wemyss}. If $A=\End(\oplus_{i \in Q_0} T_i)$ where the $T_v$ are the Cohen-Macaulays, a mutation
of $v$ will correspond to construction a new algebra $A'=\End(\oplus_{i \in Q_0\setminus \{v\}} T_i \oplus T_v')$
where $T_v'$ is the kernel of the map 
\[
b_1+b_2: T_{t(b_1)}\oplus T_{t(b_2)}\to T_v : (x,y) \mapsto x+y \text{  (note that $T_{t(b_i)}\subset T_{v}$)} 
\]
If $T_{t(b_1)}=T_{(r_1,\dots,r_k)}$ and $T_{t(b_2)}=T_{(s_1,\dots,s_k)}$ then the kernel can be identified
with the overlap of the two modules as submodules of $T_v$, so $T'_v = T_{(\max(r_1,s_1),\dots,\max(r_k,s_k))}$.
This second procedure gives the same result for the dimer model as the first, note however
that the new vertex $v'$ is not necessarily embedded on exactly the same spot as the original one. 

We have the following fact:
\begin{theorem}[Iyama-Wemyss]
The mutation of an MMA is again an MMA for the same singularity.
\end{theorem}
From this and theorems \ref{MMAisNCCR} and \ref{consistency}
we can deduce that the procedure of toric mutation turns consistent dimer models into consistent dimer models.

The analogy of this with the commutative case is remarkable. 
A toric crepant resolution of a 3-dimensional toric Gorenstein singularity corresponds to a subdivision of its polygon in elementary triangles.
Given one such subdivision one can construct a new one by looking at a quadrangle consisting of 2 two elementary triangles and
switching the diagonal so see this quadrangle as a union of 2 different triangles (for a picture see section \ref{flop}). This is called a toric flop \cite{flipflop} and
every two toric crepant resolutions can be transformed into each other by a sequence of flops.

In the noncommutative situation we have the notion of a dimer model and again there is a procedure that can transform
one dimer model into another: the toric mutation. The question now is very similar
\begin{question}\label{mutation}
Can any 2 dimer models for the same singularity be transformed into each other by a sequence of toric mutations?
\end{question}

Using our algorithm we can already answer this question affirmatively for the singularities from reflexive polygons.
\begin{theorem}
$2$ dimer models for a singularity coming from a reflexive polygon can transformed into each other by a sequence of toric mutations.
\end{theorem}
\begin{proof}
The algorithm gives us a finite list of dimer models for which we can check the statement manually.
\end{proof}

In the two following subsections we are going to explore this question for
another very special type of singularities: quotients of the conifold singularity. In this situation we will also be able to
prove a positive result.

\subsection{Quotient singularities}\label{QS}

Suppose $R=R_\sigma$ with $\sigma=\cone{v_1,\dots, v_k}\in N$.
We will identify $N$ with $\Z^n$ as column vectors and $M$ as row vectors.
Let $U$ be a matrix with integer coefficients which 
has a nonzero determinant. If all $Uv_i$ are primitive vectors, 
then we can define a new singularity with a cone $U\sigma=\cone{Uv_1,\dots ,Uv_k}$.

This new singularity is isomorphic to a quotient singularity of $\Spec R$.
Let $MU$ be the sub-lattice of $M$ of all $mU$ for which $m \in M$.
The group $G=\{\rho \in \Hom(M,\C^*)| \forall m\in MU:\rho(m)=1\}$ has an action
on $R\subset \C[M]$ by setting $(x^m)^{\rho}=\rho(m)x^m$.

From the construction it is clear that the ring of invariant functions $R^G$ 
is isomorphic to the ring $R_{U\sigma}$:
\[
 R^G = \Span_\C \{x \in MU| \<x,v_i\>\ge 0\} \cong \Span_\C \{x \in M| \<x,Uv_i\>\ge 0\}.
\]

Choose representatives in $\{m_1,\dots, m_g\}\subset M$ for the elements in the quotient group $M/MU$.
If $T_b$ is a graded reflexive $R$-submodule of $\C[M]$, then it decomposes into 
graded reflexive $R^G$-submodules 
\se{
T_b &= \Span_\C \T_b = \bigoplus_{j} \Span_\C \T_b\cap m_j+MU\\ 
&\cong \bigoplus_{j} \Span_\C \{x\in M | m_j+xU \in \T_b\}\\
&= \bigoplus_{j} T_{b,j}
}
where $T_{b,j}$ is the $R_{U\sigma}$-module $\Span_\C \{x\in M | \<x,Uv_i\> \ge b_i-\<m_j,v_i\>\}$.

\begin{lemma}\label{cover}
Let $T$ be a direct sum of graded rank one reflexive modules.
\begin{itemize}
 \item $\End_{R^G} T \cong \End_R T * G$
 \item If $\End_R T$ is an NCCR for $R$ then $\End_{R^G} T$ is an NCCR for $R^G$.
 \item If $\tilde Q \subset M\otimes \R$ is the universal cover of the embedded quiver for $\End_R T$, then
$\tilde Q\cdot U^{-1}$ will be the universal cover of the embedded quiver for $\End_{R^G}T$. 
\end{itemize}
\end{lemma}
\begin{proof}
All $T_b$ are monomial submodules of $\C[M]$ have an action of $G$ and therefore $\End_R T$ also has a $G$-action
and we can make the skew group ring $\End_R T * G$.
The orthogonal idempotents of $\End_R T * G$ are products of orthogonal idempotents in $\End_R T$ and 
orthogonal idempotents in $\C G$. The former are parametrized by the direct summands of $T=\oplus_b T_b$ while the 
latter are parametrized by the $m_i$:
\[
 e_b = \Id_{T_b} \text{ and } e_i =e_{m_i}:= \frac1{|G|}\sum_{\rho \in G} \rho(m_i) \rho.
\]
We have that $x e_{i}=e_{i-k}x$ if $x \in T_{b,k}$, so if we calculate the direct summands of $\End_R T * G$ we get
\[ e_ie_b\End_R T e_c e_j = e_i \Hom_R(T_b,T_c)e_j= e_i\bigoplus_{k}T_{c-b,k}e_j= e_iT_{c-b,j-i}
\]
which shows that $\End_R T * G\cong \End_{R^G} T$.
$\End_R T * G$ is homologically homogeneous if and only if $\End_R T$ is because the global dimension of
$\End_R T * G$ and $\End_R T$ is the same (see \cite{extgroups}).
The third statement follows from the easy to check fact that the map $\kappa$ for the new singularity is the old $\kappa$ multiplied with $U^{-1}$.
\end{proof}

\subsection{Parallellograms}
Here we will briefly examine the case where the polygon is a parallelogram. 
\[
\sigma = \cone{(0,0,1),(a,b,1),(c,d,1),(a+c,b+d,1)} 
\]
In this case the singularity can be seen
as a quotient of the conifold singularity (for which the polygon is the unit square $\cone{(0,0,1),(1,0,1),(0,1,1),(1,1,1)}$).

The conifold has a unique toric NCCR which corresponds to a torus tiled by $2$ squares.
\[
\End_R\left(T_{(0,0,0,0)}\oplus T_{(0,0,0,1)}\right)
\]
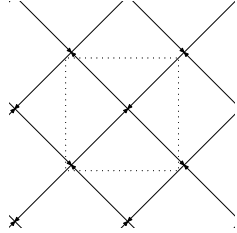
\begin{center}
\resizebox{3cm}{!}{\includegraphics{coni1.mps}}
\end{center}

By lemma \ref{cover} this quotient singularity then has an NCCR which corresponds to tiling the torus by $2 |ac-bd|$ squares,
but in general this is not the only NCCR.
 
E.g. for $\cone{(0,0,1),(2,1,1),(-1,2,1),(1,3,1)}$ we get 5 NCCRs:
\vspace{.4cm}
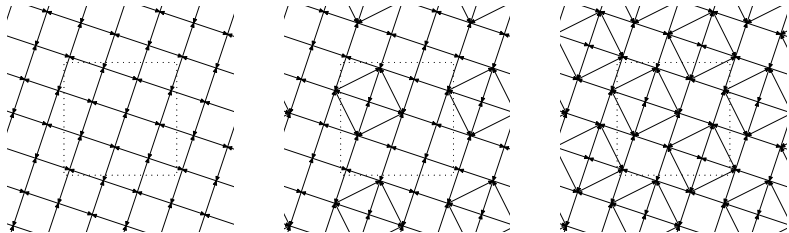
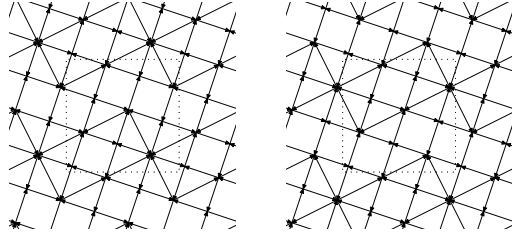
\begin{center}
\resizebox{3cm}{!}{\includegraphics{conicover5.mps}}
\hspace{.5cm}
\resizebox{3cm}{!}{\includegraphics{conicover4.mps}}
\hspace{.5cm}
\resizebox{3cm}{!}{\includegraphics{conicover3.mps}}
\hspace{.5cm}\\
\vspace{.3cm}
\resizebox{3cm}{!}{\includegraphics{conicover2.mps}}
\hspace{.5cm}
\resizebox{3cm}{!}{\includegraphics{conicover1.mps}}
\end{center}
\vspace{.4cm}
These dimers all look like the first one with
some additional \emph{diagonal arrows} that go around in curves, ignoring orientation issues.
One can solve the orientation issue by reversing the directions of the diagonal arrows pointing in the ``North-South'' directions
while keeping the orientation of the diagonal arrows in the ``East-West'' directions.
If one looks at the homology of these curves with coefficients in $\Z$, one sees
that the total homology is zero. This is not a coincidence. 

\begin{theorem}
Consider $R_\sigma$ with
$\sigma = \cone{(0,0,1),(a,b,1),(a+c,b+d,1),(c,d,1)}$.
\begin{itemize}
 \item For each dimer model of this singularity,
construct an unoriented graph by forgetting the orientations of the arrows. This graph consists
of the tiling of the torus by $2 |ac-bd|$ squares together with some additional non-intersecting diagonals. 
 \item Divide the diagonal arrows in two classes $D_+,D_-$ such that the directions of the diagonals
in each class are parallel. Then the homology class of
\[
 \sum_{a \in D_+}a - \sum_{a \in D_-}a  
\]
on the torus is zero.
\item 
All dimer models of $R$ can be torically mutated to the dimer model which tiles 
the torus by $2 |ac-bd|$ squares (i.e. the one without diagonals).
\end{itemize}
\end{theorem}
\begin{proof}
One can compute that 
\[
 \begin{pmatrix}a&b&0\\c&d&0\\0&0&1\end{pmatrix}\kappa(b_1,b_2,b_3,b_4)= 
\begin{pmatrix}
\frac{-b_1+b_2+b_3-b_4}2\\
\frac{-b_1-b_2+b_3+b_4}2\\
\frac{3b_1+b_2-b_3+b_4}4
\end{pmatrix}.
\]
In the fundamental domain of the torus $[0,1)^2$ there are only $2 |ac-bd|$ points with the property $\left(\begin{smallmatrix}a&b\\c&d\end{smallmatrix}\right)\left(\begin{smallmatrix}x\\y
\end{smallmatrix}\right)\in \frac 12\Z^2$, 
so for a dimer model, every such points must have a graded rank one Cohen-Macaulay sitting on it.
There are $4$ classes of graded rank one reflexives according to the equivalence class of $4z \mod 4$.
Note that the class to which a graded rank one reflexives belongs depends on the $(x,y)$ of its point in $[0,1)^2$ because
$2ax+2by+4z= 2b_1+2b_2=0 \mod 2$.

The directions of the arrows in the NCCR will correspond to 
$\kappa(b_1,\dots,b_4)$ 
where by corollary \ref{arrows} all $b_i=0$ except for one or two consecutive $b_i$ that are $1$. 

So in total there are eight possible directions:
The straight ones (only one $b_i\ne 0$) and the diagonal ones (2 $b_i\ne 0$).
The diagonal ones we divide in 2 classes: $D_+$ ($b_1,b_2=1$ or $b_3,b_4=1$) and
$D_-$ ($b_2,b_3=1$ or $b_1,b_4=1$). Note that after projection onto $\R^2$ all arrows in each class are (anti)parallel.
 
The consistent $\carr$-charge $\frac 12\sum_i \PM_i$ assigns degree $\frac 12$ to the straight arrows and degree $1$ to the diagonal ones.  
The consistency condition $\sum_{h(a)=v}(1-\carr_a)+\sum_{t(a)=v}(1-\carr_a)=2$ implies that every vertex has precisely 2 straight arrows arriving and 2 straight arrows leaving,
so ignoring orientation, the straight arrows form a square grid that tiles the torus.

The number of diagonal arrows in a vertex is either $2$ or $4$.
Now consider a diagonal arrow in $D_+$ with $h(a)=v$ and look at the
next diagonal arrow $b$ arriving in or leaving from $v$, in clockwise direction from $a$.
If this arrow is in $D_+$ then $h(a)=v$ because there are $2$ straight arrows between $a$ and $b$.
If this arrow is in $D_-$ then $t(a)=v$. This implies that $\sum_{a \in D_+}a - \sum_{a \in D_-}a \in \Z\qpol_1$
is a sum of cycles.

Next we show that these cycles are boundaries.
Consider $4$ points on the torus that form a grid square made of straight arrows. 
If the dimer model has a diagonal in that square then the
$2$ vertices of the square that are not on the diagonal have reflexives with different $z$-coordinates.
Indeed there is a path connecting these $2$ vertices consisting of $2$ straight arrows in the same direction and
a diagonal arrow. The $z$-coordinate of such a path is in $\frac 12 +\Z$. 
If there is no diagonal in the square or there is a diagonal connecting the two vertices then there is
a path between them consisting of 2 straight arrows. The $z$-coordinate of such a path is in $\Z$.

From this discussion we can conclude the following:
\begin{itemize}
\item the diagonal arrows connecting CMs with $4z=1,3$ form the boundary between the region containing all 
CMs with $4z=0$ and the region containing all 
CMs with $4z=2$. 
\item The diagonals going between CMs with $4z=0,2$ form 
the boundary between the region containing all 
CMs with $4z=1$ and the region containing all 
CMs with $4z=3$.
\end{itemize}
Because they are boundaries the total homology with coefficients in $\Z$ is zero.

If one mutates a vertex with coordinates $(x,y,z)\in [0,1)^3$ then the first two new coordinates
$x,y$ must be the same because the other possibilities are occupied. The third coordinate
must change by $\frac 12$ because if $(u,v,w)$ is the direction of a straight arrow and $(-u,-v,w')$ is
a straight arrow in the opposite direction then $w-w' = \pm \frac 12$. So mutation changes the equivalence class
of the vertex. 

Now look at the torus $\R^2/\Z^2$ of the dimer model and cut it in pieces along the diagonal arrows.
If there are diagonal arrows there must be at least 2 pieces because they form boundaries.

If there is a piece that is simply connected then we can 
perform a mutation on an internal vertex contained in a triangle 
(so it borders the boundary of the simply connected piece). After mutation this vertex lies outside
the simply connected piece, so gradually we can shrink the piece away.

If there is no simply connected piece, all pieces have the topology of a cylinder and are each bounded
by two curves with either the same homology class or the opposite homology class. There must be at least
one cylinder piece with opposite homology classes because the total homology of the boundaries is zero.
Using the same procedure as above we can shrink this piece away.
Eventually all pieces except one have shrunk away and we are left with the dimer model consisting only of squares.
\end{proof}

\section{Extensions of the algorithm}\label{higherdims}

We have seen that the algorithm we proposed worked for the three dimensional Gorenstein case.
We are now going to discuss briefly what happens in other cases.

\subsection{The 3-dimensional Non-Gorenstein case}\label{flop}

The concepts of NCCR and MMA were originally thought of mainly in the context of Gorenstein singularities.
Things work differently in the non-Gorenstein world. Already in the 2-dimensional situation, there is
no straightforward connection between minimal commutative resolutions and NCCRs. If we look at a 2-dimensional
abelian quotient singularity usually the NCCR, which is equal to the skew group ring, has a lot more vertices than
the minimal resolution has components and therefore the derived categories cannot be equivalent (see \cite{special}). 

In the 3-dimensional case section \ref{somefacts} still holds for non-Gorenstein singularities.
So the list $\CM$ produced by the algorithm will contain all graded rank $1$ Cohen-Macaulays. 
Step 2 will now give a list of all toric modifying modules with the highest
possible number of summands. However we do not know whether these modifying modules are maximal (it might be that higher
rank summands are needed to maximize them). We also do not know in general whether the endomorphism rings of these modules
have global dimension $3$, but often this is indeed what happens. 

Let's have a closer look at a small example: 
\[
 \tau =\cone{(0,0,1),(1,0,2),(1,1,2),(0,1,2)}
\]
which is the simplest example from toric geometry where flips occur (see \cite{flipflop}).
The algorithm produces one algebra:
\[
A_\tau := \End_{R_{\tau}} R_{\tau} \oplus T_{(0,1,1,1)}.
\]
The quiver of this algebra has 5 arrows: $3$ from the first to the second vertex 
and $2$ back. One can check that the global dimension of this algebra is $3$.

This algebra is also connected to the dimer model of the conifold.
Let $a_1,a_2,b_1,b_2$ be the $4$ arrows of the conifold dimer with
$h(a_i)=t(b_j)=v_1$ and $t(a_i)=h(b_j)=v_2$. 
By comparing generators and relations one can show that the endomorphism ring $A_\tau$ is isomorphic to the homotopy algebra of the embedded quiver that has the same vertices
as the conifold dimer but as arrows
\[
 a_1,a_2,b_2, b_1b_2^{-1}a_1^{-1}, b_1b_2^{-1}a_2^{-1}.
\]

It is well-known that the conifold algebra can be seen as the non-commutative algebra that governs the flop \cite{vdbflop}.
If one looks at the moduli space of $\theta$-stable $(1,1)$-dimensional representations \cite{Kingstab} of the conifold dimer
then the change from $\theta=(-1,1)$ to $\theta=(1,-1)$ will correspond to a flop in the moduli space.
\begin{center}
\begin{tabular}{ccc}
$\theta$ & maximal cones in the fan of moduli space\\
\hline
$(-1,1)$& $\cone{(0,0,1),(1,0,1),(1,1,1)},\cone{(0,0,1),(0,1,1),(1,1,1)}$
&\resizebox{!}{.5cm}{\begin{tikzpicture}
\filldraw [gray] 
(0,0) circle (2pt)  (1,0) circle (2pt)
(0,1) circle (2pt) (1,1) circle (2pt);
\draw (0,0) -- (1,0) -- (1,1) -- (0,1) -- (0,0);
\draw (0,0) -- (1,1);  
\end{tikzpicture}}
\\
$(0,0)$&  $\cone{(0,0,1),(1,0,1),(1,1,1),(0,1,1)}$&
\resizebox{!}{.5cm}{\begin{tikzpicture}
\filldraw [gray] 
(0,0) circle (2pt)  (1,0) circle (2pt)
(0,1) circle (2pt) (1,1) circle (2pt);
\draw (0,0) -- (1,0) -- (1,1) -- (0,1) -- (0,0);  
\end{tikzpicture}}
\\
$(1,-1)$& $\cone{0,0,1),(1,0,1),(0,1,1)},\cone{(0,1,1),(1,0,1),(1,1,1)}$&
\resizebox{!}{.5cm}{\begin{tikzpicture}
\filldraw [gray] 
(0,0) circle (2pt)  (1,0) circle (2pt)
(0,1) circle (2pt) (1,1) circle (2pt);
\draw (0,0) -- (1,0) -- (1,1) -- (0,1) -- (0,0);  
\draw (1,0) -- (0,1);  
\end{tikzpicture}}
\end{tabular}
\end{center}

Similarly one can ask whether the algebra $A_\tau$ will give you a description of the flip.
This is however not the case. If $\theta=(-1,1)$ one still gets a subdivision of $\tau$ in two smooth subcones, but
if $\theta=(1,-1)$ the cone is subdivided in 2 smooth cones and one singular
\begin{center}
\begin{tabular}{ccc}
$\theta$ & maximal cones in the fan of moduli space\\
\hline
$(-1,1)$& $\cone{(0,0,1),(1,0,2),(1,1,2)},\cone{(0,0,1),(0,1,2),(1,1,2)}$
&
\resizebox{!}{.5cm}{\begin{tikzpicture}
\filldraw [gray] 
(0,0) circle (2pt)  (1,0) circle (2pt)
(0,1) circle (2pt) (1,1) circle (2pt);
\draw (0,0) -- (1,0) -- (1,1) -- (0,1)-- (0,0);
\draw (0,0) -- (1,1);  
\end{tikzpicture}}
\\
$(0,0)$&  $\cone{(0,0,1),(1,0,2),(1,1,2),(0,1,2)}$
&\resizebox{!}{.5cm}{\begin{tikzpicture}
\filldraw [gray] 
(0,0) circle (2pt)  (1,0) circle (2pt)
(0,1) circle (2pt) (1,1) circle (2pt);
\draw (0,0) -- (1,0) -- (1,1) -- (0,1) -- (0,0);  
\end{tikzpicture}}
\\
$(1,-1)$& 
$\cone{(1,1,3),(1,0,2),(1,1,2)},
\cone{(1,1,3),(0,1,2),(1,1,2)}$
&\resizebox{!}{.5cm}{\begin{tikzpicture}
\filldraw [gray] 
(0,0) circle (2pt)  (1,0) circle (2pt)
(0,1) circle (2pt) (1,1) circle (2pt); (2/3,2/3) circle (2pt);
\draw (0,0) -- (1,0) -- (1,1) -- (0,1) -- (0,0);  
\draw (1,0) -- (2/3,2/3) -- (0,1);
\draw (1,1) -- (2/3,2/3);
\end{tikzpicture}}
\\
&
$\cone{(0,0,1),(1,0,2),(1,1,3),(0,1,2)}$&
\end{tabular}
\end{center}

\begin{remark}
These subdivisions are obtained in the following way. From \cite{lebruynbook} we know that
the moduli space can be covered by representation spaces of universal localizations of $A_\tau$. 
Each of the localizations is constructed by inverting
the nonzero paths $p$ of a $\theta$-semistable representation. 
The center of the universally localized algebra $A_\tau[p^{-1}]$ will be of the form $R_{\tau'}$ where
$\tau'$ is a subcone of $\tau$ and all these cones form a fan that subdivides $\tau$.
\end{remark}

\subsection{Higher dimensions}

In higher dimensions the algorithm we proposed does not work. This is because 
some of the lemmas that hold in the 3-dimensional case break down. 
\begin{itemize}
\item
Proposition \ref{no-bounded-cells} only holds in one direction: a lattice point in a bounded
cell will make the module not CM. But not all $T(b_1,\dots,b_k)$ that do not have
lattice points in bounded cells are CM
A counterexample is the Gorenstein pyramid 
\[\cone{(1,0,0,1),~(-1,0,0,1),~(0,1,0,0,1),~(0,-1,0,1),~(0,0,1,1)}\]
The reflexive module $T_{(1,1,-1,-1,-1)}$ is not CM but it has no bounded chambers.
\item
Proposition \ref{remove-one} is also false in higher dimensions
A counterexample is the Gorenstein octahedron.
\[
\cone{(\pm 1,0,0,1),~(0,\pm 1,0,0,1),~(0,0,\pm 1,1)}.
\]
Here $T_{(-1,-1,0,0,-1,0)}$ is CM but
$T_{(-1,-1,0,0,-1)}$ is not CM for the pyramid because $\{v_3,v_4\}$ form a disconnected subset of the pyramid.
\end{itemize}
We can partially solve this by using the singular homology computation to check whether a given 
$T_b$ is Cohen-Macaulay or not. However, as the second example above shows, the procedure
to construct a complete list of all graded rank 1 Cohen-Macaulays by increasing the rays in the fan
one by one is not exhaustive. This means that after adjusting the algorithm with singular homology computation, it
still generates an incomplete list. 

One could bypass this problem by generating a list of all normalized vectors, 
whose norm does not exceed a given number $N$ and then checking which of them are Cohen-Macaulay. 
As there are only a finite number of these, there will be a big enough $N$ such that one gets all Cohen-Macaulays.
The problem is that in general, we do not yet know a bound on the norms for Cohen-Macaulays.

Because we do not have an exhaustive list we do not know whether the modifying modules the algorithm generates
are maximal (or even just maximal for graded rank one modules).
However, to check whether the corresponding endomorphism ring is an NCCR on can always check whether this 
endomorphism ring has finite global dimension and apply lemma 4.2 of \cite{nccrep}.

In many interesting cases the algorithm does provide us with nice examples of NCCRs.
\begin{enumerate}
\item If the fan $\cone{v_1,\dots, v_n}$
is simplicial (the number of rays equals the dimension) then we are in the case of McKay correspondence.
Indeed we can apply subsection \ref{QS} for $R=\C[X_1,\dots, X_n]$ and $U=[v_1 \dots v_n]$.

As $R$ is smooth it is its own NCCR. The corresponding embedded quiver has $1$ vertex and $n$ loops.
Its universal cover
 $\tilde Q$ is the infinite embedded quiver with vertices $Q_0=M \subset M\otimes \R$ an arrows 
\[
a_{m,i}:[0,1]\to M\otimes \R: t \mapsto m+t(\underbrace{0,\dots, 0,1,0,\dots, 0}_{\text{1 on spot $i$}})
\]

The embedded quiver $\tilde Q\cdot U^{-1}$ will be the universal cover of embedded quiver of the NCCR.
The latter is hence a torus covered with $\det U$ $n$-cubes. 
A detailed description of this can be found in \cite{crawvelez}.

An interesting example of this is
\[
\cone{(0,0,0,1),(0,1,1,1),(1,0,1,1),(1,1,0,1)}
\]
because this is a singularity which allows no CCR: it is impossible to divide this tetrahedron into smaller tetrahedra because
it has no internal lattice points. It has a unique NCCR coming from
a tiling of the 4-torus with $2$ hypercubes. So this gives us a toric Gorenstein singularity without CCR but with NCCR.
If one looks at the space of $\theta$-stable representations for $\theta=(-1,1)$ or $\theta=(1,-1)$, the singularity 
resolves by dividing the cone in $4$ and introducing 1 extra ray in the $(1,1,1,2)$-direction. 
\item
We end with the 4-dimensional Gorenstein singularity generated by a unit cube.
\[
\left[\begin{smallmatrix}
(0,0,0,1),&(1,0,0,1),&(0,1,0,1),&(0,0,1,1),\\
(0,1,1,1),&(1,0,1,1),&(1,1,0,1),&(1,1,1,1) 
\end{smallmatrix}\right].
\]
This singularity is Gorenstein and can be seen as a 4-dimensional analogue of the conifold singularity.
It does have several NCCRs. One example is
given by the following maximal modifying module
\[
 T_0  \oplus T_{b_7+b_8} \oplus T_{b_6+b_8}  \oplus T_{b_6+b_7+2b_8} \oplus T_{b_8} \oplus T_{b_1}
\]
where $b_k$ is the $8$-vector with a $1$ on the $k^{th}$ entry and zero everywhere else.
The quiver can be projected to the $3$-dimensional torus, which is depicted below 
as a unit cube with the opposite faces identified.
The vertices of the quiver fall in $3$ classes. The vertex $(1)$ corresponding to the trivial CM is the vertex of the cube, 
three other vertices $(2,3,4)$ are in the centers of the faces of the cube and form the vertices of a octahedron.
They correspond to $T_{b_7+b_8}, T_{b_6+b_8}$ and $T_{b_6+b_7+2b_8}$. The two vertices $(5,6)$ coming from $T_{b_1}$ and $T_{b_8}$
lie on the diagonal of the cube, outside both sides of the octahedron and have coordinates $(\frac 14,\frac 14,\frac 14)$ and
$(\frac 34,\frac 34,\frac 34)$.

The quiver has $18$ arrows, which we can give an $\carr$-charge in a similar way to theorem \ref{charges}. 
\begin{itemize}
\item $2$ loops (vertex $5,6$) that go along the $X$-direction ($\carr$-charge $1$).
\item $4$ arrows along the diagonals of the faces of the cube from the octahedral vertices $2,3$ to the cube vertex ($\carr$-charge $\frac 12$).
\item $2$ arrows from the cube vertex to vertices $5,6$ that go along the diagonal of the cube ($\carr$-charge $\frac 14$).
\item $6$ arrows between vertex $5,6$ and the octahedral vertices ($\carr$-charge $14$).
\item $4$ arrows along the edges of the octahedron from vertices $2,3$ to $4$ ($\carr$-charge $12$).
\end{itemize}
The relations are now given by homotopy relations in the $3$-torus between paths with the same $\carr$-charge. 
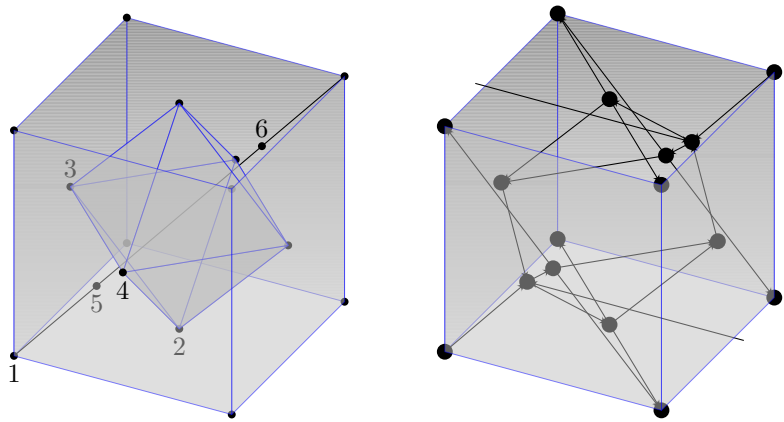
\begin{center}
\begin{tikzpicture}[y  = {(1cm,1cm)},
                    x  = {(1.93cm,-0.52cm)},
                    z  = {(0cm,2cm)},
                    scale = 1.5,
                    color = {black}]
% style of faces
\tikzset{facestyle/.style={fill=lightgray, opacity=.5,draw=blue,very thin,line join=round}}
\tikzstyle{ann} = [draw=darkgray,font=\footnotesize,inner sep=1pt]
% face "back" 
\begin{scope}[canvas is zx plane at y=1]
  \path[facestyle,shade] (0,0) rectangle (1,1);
\end{scope}
\begin{scope}[canvas is zy plane at x=0]
  \path[facestyle,shade] (0,0) rectangle (1,1);
\end{scope}
% labels

\fill(0,0,0) circle (1pt) node [below] {$1$};
\fill(1/2,1/2,0) circle (1pt) node [below] {$2$};
\fill(0,1/2,1/2) circle (1pt) node [above] {$3$};
\fill(1/2,1/2,1) circle (1pt);
\fill(1,1/2,1/2) circle (1pt);
\fill(1/2,1,1/2) circle (1pt);
\fill(1/4,1/4,1/4) circle (1pt) node [below] {$5$};
\fill(3/4,3/4,3/4) circle (1pt) node [above] {$6$};
\fill(1,0,0) circle (1pt);
\fill(0,1,0) circle (1pt);
\fill(0,0,1) circle (1pt);
\fill(1,1,0) circle (1pt);
\fill(0,1,1) circle (1pt);
\fill(1,0,1) circle (1pt);
\fill(1,1,1) circle (1pt);
\draw (0,0,0)--(1,1,1);
\path[facestyle] (1/2,1/2,0) -- (1/2,1,1/2) -- (0,1/2,1/2) -- cycle;
\path[facestyle] (1/2,1/2,1) -- (1/2,1,1/2) -- (0,1/2,1/2) -- cycle;
\path[facestyle] (1/2,1/2,0) -- (1/2,1,1/2) -- (1,1/2,1/2) -- cycle;
\path[facestyle] (1/2,1/2,1) -- (1/2,1,1/2) -- (1,1/2,1/2) -- cycle;
\path[facestyle] (1/2,1/2,0) -- (1/2,0,1/2) -- (0,1/2,1/2) -- cycle;
\path[facestyle] (1/2,1/2,1) -- (1/2,0,1/2) -- (0,1/2,1/2) -- cycle;
\path[facestyle] (1/2,1/2,0) -- (1/2,0,1/2) -- (1,1/2,1/2) -- cycle;
\path[facestyle] (1/2,1/2,1) -- (1/2,0,1/2) -- (1,1/2,1/2) -- cycle;
\begin{scope}[canvas is zy plane at x=1]
  \path[facestyle] (0,0) rectangle (1,1);
\end{scope}
% face  "right"
\begin{scope}[canvas is zx plane at y=0]
  \path[facestyle] (0,0) rectangle (1,1);
\end{scope}
\fill(1/2,0,1/2) circle (1pt) node [below] {$4$};

\end{tikzpicture}
\hspace{1cm}
\begin{tikzpicture}[y  = {(.5cm,.5cm)},
                    x  = {(.96cm,-0.26cm)},
                    z  = {(0cm,1cm)},
                    scale = 3,
                    color = {black}]
% style of faces
\tikzset{facestyle/.style={fill=lightgray, opacity=.5,draw=blue,very thin,line join=round}}
\tikzstyle{ann} = [draw=darkgray,font=\footnotesize,inner sep=1pt]
% face "back" 
\begin{scope}[canvas is zy plane at x=0]
  \path[facestyle,shade] (0,0) rectangle (1,1);
\end{scope}
% face  "left"
\begin{scope}[canvas is zx plane at y=1]
  \path[facestyle,shade] (0,0) rectangle (1,1);
\end{scope}
% face "front"

\fill(0,0,0) circle (1pt); 
\fill(1/2,1/2,0) circle (1pt); 
\fill(0,1/2,1/2) circle (1pt); 
\fill(1/2,0,1/2) circle (1pt);
\fill(1/2,1/2,1) circle (1pt);
\fill(1,1/2,1/2) circle (1pt);
\fill(1/2,1,1/2) circle (1pt);
\fill(1/4,1/4,1/4) circle (1pt); 
\fill(3/4,3/4,3/4) circle (1pt); 
\fill(1,0,0) circle (1pt);
\fill(0,1,0) circle (1pt);
\fill(0,0,1) circle (1pt);
\fill(1,1,0) circle (1pt);
\fill(0,1,1) circle (1pt);
\fill(1,0,1) circle (1pt);
\fill(1,1,1) circle (1pt);

\draw[arrows=-latex'](0,0,0)--(1/4,1/4,1/4);
\draw[arrows=-latex'](1,1,1)--(3/4,3/4,3/4);

\draw[arrows=latex'-](1,0,0)--(1/2,0,1/2);
\draw[arrows=latex'-](0,0,1)--(1/2,0,1/2);
\draw[arrows=latex'-](1,0,0)--(1/2,1/2,0);
\draw[arrows=latex'-](0,1,0)--(1/2,1/2,0);
\draw[arrows=latex'-](1,1,0)--(1/2,1,1/2);
\draw[arrows=latex'-](0,1,1)--(1/2,1,1/2);
\draw[arrows=latex'-](1,0,1)--(1/2,1/2,1);
\draw[arrows=latex'-](0,1,1)--(1/2,1/2,1);
\draw[arrows=-latex'](5/4,1/4,1/4)--(1/4,1/4,1/4);
\draw[arrows=-latex'](-1/4,3/4,3/4)--(3/4,3/4,3/4);

\draw[arrows=-latex'](1/4,1/4,1/4)--(1/2,0,1/2);
\draw[arrows=-latex'](1/4,1/4,1/4)--(1/2,1/2,0);
\draw[arrows=-latex'](3/4,3/4,3/4)--(1/2,1,1/2);
\draw[arrows=-latex'](3/4,3/4,3/4)--(1/2,1/2,1);
\draw[arrows=latex'-](1/4,1/4,1/4)--(0,1/2,1/2);
\draw[arrows=latex'-](3/4,3/4,3/4)--(1,1/2,1/2);
\draw[arrows=latex'-](0,1/2,1/2)--(1/2,1,1/2);
\draw[arrows=latex'-](1,1/2,1/2)--(1/2,0,1/2);
\draw[arrows=latex'-](0,1/2,1/2)--(1/2,1/2,1);
\draw[arrows=latex'-](1,1/2,1/2)--(1/2,1/2,0);

\begin{scope}[canvas is zy plane at x=1]
  \path[facestyle] (0,0) rectangle (1,1);
\end{scope}
% face  "right"
\begin{scope}[canvas is zx plane at y=0]
  \path[facestyle] (0,0) rectangle (1,1);
\end{scope}
\end{tikzpicture}
\end{center}
\end{enumerate}

\section{Appendix: Dimer models for Fano surfaces}

In this appendix we give a list of all reflexive polygons, the corresponding dimer models and their type.
Some dimer models are not isomorphic to their opposite (i.e. all arrows reversed). 
We will include only one of the $2$ orientations and denote 
dimer models with an asterisk. 

\begin{itemize}
\item [\resizebox{!}{1cm}{\begin{tikzpicture}
\filldraw [gray] 
(-.5,-.5) circle (2pt) (0,-.5) circle (2pt) (.5,-.5) circle (2pt)
(-.5,0) circle (2pt)  (.5,0) circle (2pt)
(-.5,.5) circle (2pt) (0,.5) circle (2pt) (.5,.5) circle (2pt);
\draw (0,0) node {3a};
\draw (-.5,-.5) -- (.5,0) -- (0,.5) -- (-.5,-.5);  
\end{tikzpicture}}]
\begin{tabular}{ccc}
\resizebox{3cm}{!}{\includegraphics{refl1.mps}}&&\\
3a&&
\end{tabular}
\item [\resizebox{!}{1cm}{\begin{tikzpicture}
\filldraw [gray] 
(-.5,-.5) circle (2pt) (0,-.5) circle (2pt) (.5,-.5) circle (2pt)
(-.5,0) circle (2pt)  (.5,0) circle (2pt)
(-.5,.5) circle (2pt) (0,.5) circle (2pt) (.5,.5) circle (2pt);
\draw (0,-.5) -- (.5,0) -- (0,.5) -- (-.5,0) -- (0,-.5);  
\draw (0,0) node {4a};
\end{tikzpicture}}]
\begin{tabular}{ccc}
\resizebox{3cm}{!}{\includegraphics{refl2.mps}}&
\resizebox{3cm}{!}{\includegraphics{refl3.mps}}&\\
4a&4c&\\
\end{tabular}
\item [\resizebox{!}{1cm}{\begin{tikzpicture}
\filldraw [gray] 
(-.5,-.5) circle (2pt) (0,-.5) circle (2pt) (.5,-.5) circle (2pt)
(-.5,0) circle (2pt)  (.5,0) circle (2pt)
(-.5,.5) circle (2pt) (0,.5) circle (2pt) (.5,.5) circle (2pt);
\draw (.5,-.5) -- (0,.5) -- (-.5,0) -- (0,-.5) -- (.5,-.5);  
\draw (0,0) node {4b};
\end{tikzpicture}}]
\begin{tabular}{ccc}
\resizebox{3cm}{!}{\includegraphics{refl4.mps}}&&\\
4b&&
\end{tabular}
\item [\resizebox{!}{1cm}{\begin{tikzpicture}
\filldraw [gray] 
(-.5,-.5) circle (2pt) (0,-.5) circle (2pt) (.5,-.5) circle (2pt)
(-.5,0) circle (2pt)  (.5,0) circle (2pt)
(-.5,.5) circle (2pt) (0,.5) circle (2pt) (.5,.5) circle (2pt);
\draw (.5,-.5) -- (0,.5) -- (-.5,-.5) -- (.5,-.5);  
\draw (0,0) node {4c};
\end{tikzpicture}}]
\begin{tabular}{ccc}
\resizebox{3cm}{!}{\includegraphics{refl5.mps}}&&\\
4a&&
\end{tabular}
\item [\resizebox{!}{1cm}{\begin{tikzpicture}
\filldraw [gray] 
(-.5,-.5) circle (2pt) (0,-.5) circle (2pt) (.5,-.5) circle (2pt)
(-.5,0) circle (2pt)  (.5,0) circle (2pt)
(-.5,.5) circle (2pt) (0,.5) circle (2pt) (.5,.5) circle (2pt);
\draw (.5,0) -- (0,.5) -- (-.5,.5) -- (-.5,0) -- (0,-.5) -- (.5,0);  
\draw (0,0) node {5a};
\end{tikzpicture}}]
\begin{tabular}{ccc}
\resizebox{3cm}{!}{\includegraphics{refl6.mps}}&
\resizebox{3cm}{!}{\includegraphics{refl7.mps}}&\\
5a&5b*&
\end{tabular}
\item [\resizebox{!}{1cm}{\begin{tikzpicture}
\filldraw [gray] 
(-.5,-.5) circle (2pt) (0,-.5) circle (2pt) (.5,-.5) circle (2pt)
(-.5,0) circle (2pt)  (.5,0) circle (2pt)
(-.5,.5) circle (2pt) (0,.5) circle (2pt) (.5,.5) circle (2pt);
\draw (.5,-.5) -- (0,.5) -- (-.5,0) -- (-.5,-.5) -- (.5,-.5);  
\draw (0,0) node {5b};
\end{tikzpicture}}]
\begin{tabular}{ccc}
\resizebox{3cm}{!}{\includegraphics{refl9.mps}}&&\\
5a&&
\end{tabular}
\item [\resizebox{!}{1cm}{\begin{tikzpicture}
\filldraw [gray] 
(-.5,-.5) circle (2pt) (0,-.5) circle (2pt) (.5,-.5) circle (2pt)
(-.5,0) circle (2pt)  (.5,0) circle (2pt)
(-.5,.5) circle (2pt) (0,.5) circle (2pt) (.5,.5) circle (2pt);
\draw (.5,0) -- (0,.5) -- (-.5,.5) -- (-.5,0) -- (0,-.5) -- (.5,-.5) -- (.5,0);  
\draw (0,0) node {6a};
\end{tikzpicture}}]
\begin{tabular}{ccc}
\resizebox{3cm}{!}{\includegraphics{refl10.mps}}&
\resizebox{3cm}{!}{\includegraphics{refl11.mps}}&
\resizebox{3cm}{!}{\includegraphics{refl12.mps}}\\
6b&6a&6c\\
\resizebox{3cm}{!}{\includegraphics{refl13.mps}}&
\resizebox{3cm}{!}{\includegraphics{refl14.mps}}&\\
6c&6d*&\\
\end{tabular}
\item [\resizebox{!}{1cm}{\begin{tikzpicture}
\filldraw [gray] 
(-.5,-.5) circle (2pt) (0,-.5) circle (2pt) (.5,-.5) circle (2pt)
(-.5,0) circle (2pt)  (.5,0) circle (2pt)
(-.5,.5) circle (2pt) (0,.5) circle (2pt) (.5,.5) circle (2pt);
\draw (.5,0) -- (0,.5) -- (-.5,.5) -- (-.5,-.5) -- (0,-.5) -- (.5,0);  
\draw (0,0) node {6b};
\end{tikzpicture}}]
\begin{tabular}{ccc}
\resizebox{3cm}{!}{\includegraphics{refl16.mps}}&
\resizebox{3cm}{!}{\includegraphics{refl17.mps}}&
\resizebox{3cm}{!}{\includegraphics{refl18.mps}}\\
6c*&6b&6a
\end{tabular}
\item [\resizebox{!}{1cm}{\begin{tikzpicture}
\filldraw [gray] 
(-.5,-.5) circle (2pt) (0,-.5) circle (2pt) (.5,-.5) circle (2pt)
(-.5,0) circle (2pt)  (.5,0) circle (2pt)
(-.5,.5) circle (2pt) (0,.5) circle (2pt) (.5,.5) circle (2pt)
(-1,-.5) circle (2pt) (-1,0) circle (2pt) (-1,.5) circle (2pt);
\draw (.5,0) -- (0,.5) -- (-1,-.5) -- (0,-.5) -- (.5,0);  
\draw (0,0) node {6c};
\end{tikzpicture}}]
\begin{tabular}{ccc}
\resizebox{3cm}{!}{\includegraphics{refl20.mps}}&
\resizebox{3cm}{!}{\includegraphics{refl21.mps}}&\\
6b&6a&
\end{tabular}
\item [\resizebox{!}{1cm}{\begin{tikzpicture}
\filldraw [gray] 
(-.5,-.5) circle (2pt) (0,-.5) circle (2pt) (.5,-.5) circle (2pt)
(-.5,0) circle (2pt)  (.5,0) circle (2pt)
(-.5,.5) circle (2pt) (0,.5) circle (2pt) (.5,.5) circle (2pt)
(-1,-.5) circle (2pt) (-1,0) circle (2pt) (-1,.5) circle (2pt);
\draw (0,.5) -- (-1,-.5) -- (.5,-.5) -- (0,.5);  
\draw (0,0) node {6d};
\end{tikzpicture}}]
\begin{tabular}{ccc}
\resizebox{3cm}{!}{\includegraphics{refl22.mps}}&
&\\
6a&&\\
\end{tabular}
\item [\resizebox{!}{1cm}{\begin{tikzpicture}
\filldraw [gray] 
(-.5,-.5) circle (2pt) (0,-.5) circle (2pt) (.5,-.5) circle (2pt)
(-.5,0) circle (2pt)  (.5,0) circle (2pt)
(-.5,.5) circle (2pt) (0,.5) circle (2pt) (.5,.5) circle (2pt);
\draw (.5,0) -- (0,.5) -- (-.5,.5) -- (-.5,-.5) -- (.5,-.5) -- (.5,0); 
\draw (0,0) node {7a};
\end{tikzpicture}}]
\begin{tabular}{ccc}
\resizebox{3cm}{!}{\includegraphics{refl23.mps}}&
\resizebox{3cm}{!}{\includegraphics{refl24.mps}}&
\resizebox{3cm}{!}{\includegraphics{refl27.mps}}\\
7b*&7a*&7a\\
\end{tabular}
\item [\resizebox{!}{1cm}{\begin{tikzpicture}
\filldraw [gray] 
(-.5,-.5) circle (2pt) (0,-.5) circle (2pt) (.5,-.5) circle (2pt)
(-.5,0) circle (2pt)  (.5,0) circle (2pt)
(-.5,.5) circle (2pt) (0,.5) circle (2pt) (.5,.5) circle (2pt)
(-1,-.5) circle (2pt) (-1,0) circle (2pt) (-1,.5) circle (2pt);
\draw (.5,0) -- (0,.5) -- (-1,-.5) -- (.5,-.5) --(.5,0); 
\draw (0,0) node {7b};
\end{tikzpicture}}]
\begin{tabular}{ccc}
\resizebox{3cm}{!}{\includegraphics{refl28.mps}}&
&\\
7a\\
\end{tabular}
\item [\resizebox{!}{1cm}{\begin{tikzpicture}
\filldraw [gray] 
(-.5,-.5) circle (2pt) (0,-.5) circle (2pt) (.5,-.5) circle (2pt)
(-.5,0) circle (2pt)  (.5,0) circle (2pt)
(-.5,.5) circle (2pt) (0,.5) circle (2pt) (.5,.5) circle (2pt);
\draw (.5,.5) -- (-.5,.5) -- (-.5,-.5) -- (.5,-.5) -- (.5,.5);  
\draw (0,0) node {8a};
\end{tikzpicture}}]
\begin{tabular}{ccc}
\resizebox{3cm}{!}{\includegraphics{refl29.mps}}&
\resizebox{3cm}{!}{\includegraphics{refl30.mps}}&
\resizebox{3cm}{!}{\includegraphics{refl31.mps}}\\
8a&8b&8c\\
\resizebox{3cm}{!}{\includegraphics{refl32.mps}}&&
\\
8a&\\
\end{tabular}
\item [\resizebox{!}{1cm}{\begin{tikzpicture}
\filldraw [gray] 
(-.5,-.5) circle (2pt) (0,-.5) circle (2pt) (.5,-.5) circle (2pt)
(-.5,0) circle (2pt)  (.5,0) circle (2pt)
(-.5,.5) circle (2pt) (0,.5) circle (2pt) (.5,.5) circle (2pt)
(1,-.5) circle (2pt) (1,0) circle (2pt) (1,.5) circle (2pt);
\draw (0,.5) -- (-.5,.5) -- (-.5,-.5) -- (1,-.5) -- (0,.5);  
\draw (0,0) node {8b};
\end{tikzpicture}}]
\begin{tabular}{ccc}
\resizebox{3cm}{!}{\includegraphics{refl33.mps}}&
\resizebox{3cm}{!}{\includegraphics{refl34.mps}}&\\
8b&8a&\\
\end{tabular}
\item [\resizebox{!}{1cm}{\begin{tikzpicture}
\filldraw [gray] 
(-.5,-.5) circle (2pt) (0,-.5) circle (2pt) (.5,-.5) circle (2pt)
(-.5,0) circle (2pt)  (.5,0) circle (2pt)
(-.5,.5) circle (2pt) (0,.5) circle (2pt) (.5,.5) circle (2pt)
(1,-.5) circle (2pt) (1,0) circle (2pt) (1,.5) circle (2pt)
(-1,-.5) circle (2pt) (-1,0) circle (2pt) (-1,.5) circle (2pt);
\draw (0,.5) -- (1,-.5) -- (-1,-.5) -- (0,.5);  
\draw (0,0) node {8c};
\end{tikzpicture}}]
\begin{tabular}{ccc}
\resizebox{3cm}{!}{\includegraphics{refl35.mps}}&&\\
8a&&\\
\end{tabular}
\item [\resizebox{!}{1cm}{\begin{tikzpicture}
\filldraw [gray] 
(-.5,-.5) circle (2pt) (0,-.5) circle (2pt) (.5,-.5) circle (2pt)
(-.5,0) circle (2pt)  (.5,0) circle (2pt)
(-.5,.5) circle (2pt) (0,.5) circle (2pt) (.5,.5) circle (2pt)
(1,-.5) circle (2pt) (1,0) circle (2pt) (1,.5) circle (2pt)
(-.5,1) circle (2pt) (0,1) circle (2pt) (.5,1) circle (2pt)
(1,1) circle (2pt);
\draw (-.5,1) -- (-.5,-.5) -- (1,-.5) -- (-.5,1);  
\draw (0,0) node {9a};
\end{tikzpicture}}]
\begin{tabular}{ccc}
\resizebox{3cm}{!}{\includegraphics{refl36.mps}}&&\\
9a&&\\
\end{tabular}
\end{itemize}

\bibliographystyle{amsplain}
\def\cprime{$'$}
\providecommand{\bysame}{\leavevmode\hbox to3em{\hrulefill}\thinspace}
\providecommand{\MR}{\relax\ifhmode\unskip\space\fi MR }
\providecommand{\MRhref}[2]{%
  \href{http://www.ams.org/mathscinet-getitem?mr=#1}{#2}
}
\providecommand{\href}[2]{#2}

\end{document}